\documentclass[journal,twoside,web]{ieeecolor}

\usepackage{generic}
\usepackage{amsmath,mathrsfs,amsfonts,amssymb,graphicx,epsfig}
 \usepackage{amsthm}
\usepackage{subcaption}
\usepackage{color,multirow,rotating}
\usepackage{algorithm,algpseudocode,algorithmicx}
\usepackage{cite,url,framed,bm,balance,nicematrix}
\usepackage{stmaryrd,hyperref}

\newtheorem{theorem}{Theorem}
\newtheorem{corollary}[theorem]{Corollary}

\newtheorem{lemma}{Lemma}
\newtheorem{proposition}{Proposition}
\newtheorem{remark}{Remark}

\newtheorem{example}{Example}

\usepackage{tikz}
\usetikzlibrary{calc,trees,positioning,arrows,chains,shapes.geometric,decorations.pathreplacing,decorations.pathmorphing,shapes, matrix,shapes.symbols}

\newcommand{\differential}{{\rm{d}}}

\renewcommand{\det}{{\mathrm{det}}}
\newcommand{\tr}{{\mathrm{trace}}}

\newcommand{\dist}{{\mathrm{dist}}}

\makeatletter
\newcommand{\pushright}[1]{\ifmeasuring@#1\else\omit\hfill$\displaystyle#1$\fi\ignorespaces}
\makeatother

\makeatletter
\newcommand{\dotminus}{\mathbin{\text{\@dotminus}}}

\newcommand{\@dotminus}{%
  \ooalign{\hidewidth\raise1ex\hbox{.}\hidewidth\cr$\m@th-$\cr}%
}

\allowdisplaybreaks

\def\BibTeX{{\rm B\kern-.05em{\sc i\kern-.025em b}\kern-.08em
    T\kern-.1667em\lower.7ex\hbox{E}\kern-.125emX}}
\begin{document}
\title{Markov Kernels, Distances and Optimal Control: A Parable of Linear Quadratic Non-Gaussian Distribution Steering}
\author{Alexis M.H. Teter, Wenqing Wang, Sachin Shivakumar, Abhishek Halder,~\IEEEmembership{Senior Member,~IEEE}
\thanks{Alexis M.H. Teter is with the Department of Applied Mathematics, University of California, Santa Cruz, CA 95064, USA, {\tt\small{amteter@ucsc.edu}}.}%
\thanks{Wenqing Wang and Abhishek Halder are with the Department of Aerospace Engineering, Iowa State University, Ames, IA 50011, USA, {\tt\small{\{wqwang,ahalder\}@iastate.edu}}.}%
\thanks{Sachin Shivakumar is with the Theoretical Division, Los Alamos National Laboratory, Los Alamos, NM 87545, USA, {\tt\small{sshivakumar@lanl.gov}}.}%
}

\maketitle

\begin{abstract}
For a controllable linear time-varying (LTV) pair $(\bm{A}_t,\bm{B}_t)$ and $\bm{Q}_{t}$ positive semidefinite, we derive the Markov kernel $\kappa$ for the It\^{o} diffusion $$\differential\bm{x}_{t}=\bm{A}_{t}\bm{x}_t \differential t + \sqrt{2}\bm{B}_{t}\differential\bm{w}_{t}$$ with an accompanying killing of probability mass at rate $\frac{1}{2}\bm{x}^{\top}\bm{Q}_{t}\bm{x}$. This Markov kernel is the Green's function for the linear reaction-advection-diffusion partial differential equation 
$$\partial_{t}\kappa = -\langle\nabla_{\bm{x}},\kappa\bm{A}_{t}\bm{x}\rangle + \langle\bm{B}_{t}\bm{B}_{t}^{\top},\nabla_{\bm{x}}^{2}\kappa\rangle - \frac{1}{2}\bm{x}^{\top}\bm{Q}_t\bm{x}\kappa.$$
Our result generalizes the recently derived kernel for the special case  $\left(\bm{A}_t,\bm{B}_t\right)=\left(\bm{0},\bm{I}\right)$, and depends on the solution of an associated Riccati matrix ODE. A consequence of this result is that the linear quadratic non-Gaussian Schr\"{o}dinger bridge is exactly solvable. This means that the problem of steering a controlled LTV diffusion from a given non-Gaussian distribution to another over a fixed deadline while minimizing an expected quadratic cost can be solved using dynamic Sinkhorn recursions performed with the derived kernel. The endpoint non-Gaussian distributions are only required to have finite second moments, and are arbitrary otherwise. 

Our derivation for the $\left(\bm{A}_t,\bm{B}_t\right)$-parametrized kernel pursues a new idea that relies on finding a state-time dependent distance-like functional given by the solution of a deterministic optimal control problem. This technique breaks away from existing methods, such as generalizing Hermite polynomials or Weyl calculus, which have seen limited success in the reaction-diffusion context. Our technique uncovers a new connection between Markov kernels, distances, and optimal control. This connection is of interest beyond its immediate application in solving the linear quadratic Schr\"{o}dinger bridge problem. 
\end{abstract}

\begin{IEEEkeywords}
Markov kernel, stochastic optimal control, non-Gaussian distribution, Green's function, Schr\"{o}dinger bridge.
\end{IEEEkeywords}


\section{Introduction}\label{sec:introduction}
Markov kernels 
$$\kappa(t_{0},\bm{x},t,\bm{y}), \quad 0\leq t_{0} \leq t < \infty, \quad \bm{x},\bm{y}\in\mathbb{R}^{n},$$ 
play a central role in the analysis \cite[Ch. 1]{stroock2008partial} and control \cite{caluya2021wasserstein,chen2021stochastic,teter2023contraction,teter2024schr}, \cite[Ch. V-VI]{fleming2012deterministic} of Markov diffusion processes. Often but not always, they can be interpreted as \emph{transition probability}, which are measurable maps sending \emph{Borel probability measures}\footnote{endowed with the topology of weak convergence} supported on subsets of $\mathbb{R}^{n}$ to itself. A well-known example is the (Euclidean) \emph{heat kernel} \cite[p. 44-47]{evans2022partial}
\begin{align}
\kappa_{\mathrm{Heat}}(t_{0},\bm{x},t,\bm{y}) := \frac{1}{\left(4 \pi\left(t-t_0\right)\right)^{n / 2}}
\exp \left(-\frac{|\boldsymbol{x}-\boldsymbol{y}|^2}{4\left(t-t_0\right)}\right),
\label{HeatKernel} 
\end{align}
which is the transition probability for the It\^{o} process
\begin{align}
\differential\bm{x}_{t} = \sqrt{2}\:\differential\bm{w}_{t}, \quad \bm{w}_{t}\in\mathbb{R}^{n},
\label{BrownianSDE}
\end{align}
where $\bm{w}_{t}$ is the standard Wiener process, and $\vert\cdot\vert$ denotes the Euclidean norm. 

A more general example of interest in systems-control is the Markov kernel
\begin{align}
&\kappa_{\mathrm{Linear}}\left(t_0, \boldsymbol{x}, t, \boldsymbol{y}\right)=\left(4 \pi\left(t-t_0\right)\right)^{-n/2} \operatorname{det}\left(\boldsymbol{\Gamma}_{tt_0}\right)^{-1 / 2} \nonumber\\
& \qquad\qquad\exp \left(-\frac{\left(\boldsymbol{\Phi}_{tt_0} \boldsymbol{x}-\boldsymbol{y}\right)^{\top} \boldsymbol{\Gamma}_{t t_0}^{-1}\left(\boldsymbol{\Phi}_{tt_0} \boldsymbol{x}-\boldsymbol{y}\right)}{4\left(t-t_0\right)}\right),
\label{LinearKernel}    
\end{align}
which is the transition probability for the It\^{o} process
\begin{align}
\differential\bm{x}_{t} = \bm{A}_{t}\bm{x}_t\:\differential t + \sqrt{2}\:\bm{B}_{t}\differential\bm{w}_{t},\quad \bm{w}_{t}\in\mathbb{R}^{m},
\label{LinearSDE}
\end{align}
where $\left(\bm{A}_t,\bm{B}_t\right)$ is a uniformly controllable matrix-valued trajectory pair in $\left(\mathbb{R}^{n\times n},\mathbb{R}^{n\times m}\right)$ that is bounded and continuous in $t\in[t_0,\infty)$, and the state transition matrix $\bm{\Phi}_{t \tau}:=\bm{\Phi}(t, \tau) \forall t_0 \leq \tau \leq t$. Here, uniformly controllable means positive definiteness of the associated finite horizon controllability Gramian $\boldsymbol{\Gamma}_{t t_0}$, i.e.,
$$\boldsymbol{\Gamma}_{t t_0}:=\!\!\int_{t_0}^{t}\! \boldsymbol{\Phi}_{t\tau} \boldsymbol{B}_\tau \boldsymbol{B}_\tau^{\top} \boldsymbol{\Phi}_{t \tau}^{\top} \mathrm{d} \tau \succ \mathbf{0},\quad 0\leq t_{0}\leq t<\infty.$$
As expected, \eqref{LinearKernel} reduces to \eqref{HeatKernel} for  $\left(\bm{A}_t,\bm{B}_t\right)=\left(\bm{0},\bm{I}\right)$.

Both \eqref{HeatKernel} and \eqref{LinearKernel} are instances of $\kappa$ that are transition probabilities, and satisfy $\kappa\geq 0$, $\int_{\mathbb{R}^{n}}\kappa\differential\bm{y}=1$. They solve \emph{Kolmogorov's forward partial differential equation} (PDE) initial value problem with Dirac delta initial condition:
\begin{align}
\partial_{t}\kappa = L\kappa, \quad \kappa(t_{0},\bm{x},t_{0},\bm{y})=\delta(\bm{x}-\bm{y}),
\label{KolForwardPDEIVP}   
\end{align}
where $L$ is an \emph{advection-diffusion} spatial operator induced by the drift and diffusion coefficients of the underlying It\^{o} process. For example, for the It\^{o} process \eqref{BrownianSDE}, $L\kappa\equiv\Delta_{\bm{x}}\kappa$ (standard Laplacian). For \eqref{LinearSDE}, $L\kappa\equiv -\langle\nabla_{\bm{x}},\kappa\bm{A}_{t}\bm{x}\rangle + \bm{B}_{t}\bm{B}_{t}^{\top}\Delta_{\bm{x}}\kappa$. 

More generally, for the It\^{o} process $$\differential\bm{x}_{t}=\bm{f}(t,\bm{x}_{t})\differential t + \bm{g}(t,\bm{x}_{t})\differential\bm{w}_{t},$$ with Lipschitz $\bm{f}$, uniformly lower bounded $\bm{G}:=\bm{gg}^{\top}$, and $\vert \bm{f}\vert + \vert\bm{g}\vert \leq c\left(1 + \vert \bm{x}\vert\right)$ uniformly in $t$ for some constant $c>0$, we have $L\kappa\equiv \langle\nabla_{\bm{x}},\kappa\bm{f}\rangle + \Delta_{\bm{G}}\kappa$ where the weighted Laplacian $\Delta_{\bm{G}}:=\sum_{i,j}\partial^{2}_{x_{i}x_{j}}\left(\kappa G_{ij}\right)$. At this level of generality, closed-form formula for the transition probability $\kappa$ such as \eqref{HeatKernel} or \eqref{LinearKernel} are not available.

A typical situation where the Markov kernel $\kappa$ is \emph{not a transition probability} is when the underlying It\^{o} process, in addition to drift and diffusion, also allows the creation or killing of probability mass at a rate $q(\bm{x}_t)$ for some bounded measurable $q$. In such cases, $\kappa$ is a measurable map that sends \emph{nonnegative Borel measures}\footnote{still endowed with the topology of weak convergence} supported on subsets of $\mathbb{R}^{n}$ to itself. Then, \eqref{KolForwardPDEIVP} gets replaced with
\begin{align}
\partial_{t}\kappa = \left(L-q\right)\kappa, \quad \kappa(t_{0},\bm{x},t_{0},\bm{y})=\delta(\bm{x}-\bm{y}),
\label{ReactionAdvectionDiffusionPDEIVP}    
\end{align}
where $L$ is an advection-diffusion operator as before ($q=0$ case) and $L-q$ becomes a \emph{reaction-advection-diffusion} operator. We say that $q$ is the reaction rate.

For both \eqref{KolForwardPDEIVP} and \eqref{ReactionAdvectionDiffusionPDEIVP}, the Markov kernel $\kappa$ can be seen as the \emph{Green's function} of the associated linear PDE initial value problem. Thus, a closed-form handle on $\kappa$ helps solve the associated prediction problem in the sense if the initial state $\bm{x}_{0}\sim\mu_{0}$ (symbol $\sim$ denotes ``follows the statistics of") then $\bm{x}_{t}\sim \int_{\mathbb{R}^{n}} \kappa(t_0,\bm{x},t,\bm{y})\differential\mu_{0}(\bm{y})$.

Compared to \eqref{KolForwardPDEIVP}, explicit formula for $\kappa$ in \eqref{ReactionAdvectionDiffusionPDEIVP} are less known. Recently, a closed-form formula for $\kappa$ in the case $L\equiv\Delta_{\bm{x}}$ and convex quadratic $q(\bm{x}_{t})$ was found using Hermite polynomials \cite{teter2024schr} and Weyl calculus \cite{teter2024weyl}. Motivation behind these studies came from the \emph{Schr\"{o}dinger bridge problem} (SBP) with quadratic state cost $q(\bm{x}_{t}):=\frac{1}{2}\bm{x}_t^{\top}\bm{Q}\bm{x}_t$, $\bm{Q}\succeq\bm{0}$, which are stochastic optimal control problems of the form:
\begin{subequations}
\begin{align}
&\underset{\left(\mu^{\bm{u}},\bm{u}\right)}{\inf}\int_{\mathbb{R}^{n}}\int_{t_{0}}^{t_{1}}\bigg\{\frac{1}{2}\vert\bm{u}\vert^{2} + \frac{1}{2}\left(\bm{x}_t^{\bm{u}}\right)^{\top}\bm{Q}\bm{x}_t^{\bm{u}}\bigg\}\differential t \:\differential\mu^{\bm{u}}(\bm{x}_{t}^{\bm{u}})\label{ClassicalSBPstatecostObj}\\
&\text{subject to}\quad \differential\bm{x}_{t}^{\bm{u}} = \bm{u}_{t}\left(t,\bm{x}_t^{\bm{u}}\right)\:\differential t + \sqrt{2}\:\differential\bm{w}_{t},\label{ClassicalSBPstatecostSDE}\\
&\qquad\qquad\quad\bm{x}_{t}^{\bm{u}}(t=t_0) \sim \mu_{0}, \quad \bm{x}_{t}^{\bm{u}}(t=t_1) \sim \mu_{1}, \label{ClassicalSBPstatecostEndpointConstr}   
\end{align}
\label{ClassicalSBPstatecost}
\end{subequations}
where the deadline $[t_{0},t_{1}]$, the state cost weight matrix $\bm{Q}\succeq\bm{0}$, and the endpoint statistics $\mu_0,\mu_1$ are given as problem data. 

Problem \eqref{ClassicalSBPstatecost} has the interpretation of linear quadratic (LQ) optimal control synthesis for steering \emph{non-Gaussian state statistics} over a given finite horizon. The existence-uniqueness of solution for \eqref{ClassicalSBPstatecost} is guaranteed provided $\mu_0,\mu_1$ have finite second moments. Solution of \eqref{ClassicalSBPstatecost} for \emph{Gaussian} $\mu_0,\mu_1$ was detailed in \cite[Sec. III]{7170905} in a more general setting with \eqref{ClassicalSBPstatecostSDE} replaced by the controlled linear time-varying (LTV) dynamics: 
\begin{align}
\differential\bm{x}_{t}^{\bm{u}} = \left(\bm{A}_t\bm{x}_t^{\bm{u}} +\bm{B}_t\bm{u}_{t}\right)\:\differential t + \sqrt{2}\bm{B}_{t}\:\differential\bm{w}_{t},
\label{LinearSBPSDE}
\end{align}
wherein as before, $\left(\bm{A}_t,\bm{B}_t\right)$ is controllable matrix-valued trajectory pair in $\left(\mathbb{R}^{n\times n},\mathbb{R}^{n\times m}\right)$ that is bounded and continuous in $t\in[t_0,t_{1}]$.

Writing the conditions of optimality for \eqref{ClassicalSBPstatecost} followed by certain change-of-variables, it can be shown \cite[Sec. 3]{teter2024schr} that solving problem \eqref{ClassicalSBPstatecost} leads to computing the ``propagator" a.k.a. the action of the Green's function 
$$\int_{\mathbb{R}^{n}} \kappa(t_0,\bm{x},t,\bm{y})\widehat{\varphi}_{0}(\bm{y})\differential\bm{y},$$ 
where $\kappa$ solves \eqref{ReactionAdvectionDiffusionPDEIVP}, and $\widehat{\varphi}_{0}$ is a suitable measurable function. In other words, \emph{the state cost-to-go manifests as a reaction rate in the PDE for the Markov kernel}. Then, knowing a closed-form formula for $\kappa$ facilitates the solution of \eqref{ClassicalSBPstatecost} with \emph{generic non-Gaussian} $\mu_0,\mu_1$ having finite second moments. This is what was accomplished in \cite{teter2024schr,teter2024weyl}.

A natural question is whether such a closed-form formula for $\kappa$ can be derived when \eqref{ClassicalSBPstatecostSDE} is replaced by \eqref{LinearSBPSDE}. Finding such $\kappa$ would enable solving \emph{LQ SBPs with generic non-Gaussian} $\mu_0,\mu_1$ having finite second moments. From a probabilistic point of view, this $\kappa$ is the Markov kernel of the It\^{o} process \eqref{LinearSDE} with quadratic creation or killing of probability mass with rate $q(\bm{x}_{t}):=\frac{1}{2}\bm{x}_t^{\top}\bm{Q}\bm{x}_t$, $\bm{Q}\succeq\bm{0}$. In this work, we derive this Markov kernel in the more general setting of time-varying weight matrices (see Assumptions \textbf{A1}-\textbf{A2} in Sec. \ref{sec:MainSection}). 

\subsubsection*{Contributions}
This work makes two concrete contributions.

\begin{itemize}
\item The first contribution is the solution of a specific problem. We deduce the Markov kernel for the It\^{o} diffusion \eqref{LinearSDE} with killing rate $\frac{1}{2}\bm{x}^{\top}\bm{Q}_{t}\bm{x}$ for $\bm{Q}_{t}\succeq\bm{0}$. We explain how integral transforms defined by this kernel help in solving the generic LQ SBP.

\vspace*{0.1in}

\item The second contribution is methodological. To derive the aforesaid kernel, we propose a new technique that involves identifying a deterministic optimal control problem from the It\^{o} diffusion, solving the same to find a distance function, and then to identify the Markov kernel in a structure defined by the same. We provide computational details to demonstrate that the proposed technique systematically recovers Markov kernels of interest: old and new.  
\end{itemize}

\subsubsection*{Related works}
In Table \ref{table:relatedworksLQSBP}, we contrast the technical contribution of this work vis-\`{a}-vis the related works in the literature.

Of particular relevance are \cite{chen2015optimal,chenPartIII,teter2024schr,teter2024weyl}, all of which consider the quadratic state cost in the cost-to-go. The developments in \cite{chen2015optimal,chenPartIII} focused on $\mu_0,\mu_1$ Gaussian, and that case did not require the kernel, thanks to linear dynamics preserving the Gaussianity. On the other hand, the works \cite{teter2024schr,teter2024weyl} derived the kernel only for $\left(\bm{A}_t,\bm{B}_t\right)=(\bm{0},\bm{I})$ and constant $\bm{Q}\succeq 0$.

Neither the Hermite polynomials in \cite{teter2024schr} nor the Weyl calculus computation in \cite{teter2024weyl} generalize in any obvious way for deriving the $\kappa$ in the LQ setting of our interest. In this work, we chart a new path motivated by a basic observation on the structure of the Markov kernels, distances and certain deterministic optimal control problems.

\begin{table}[t]
\centering
\begin{tabular}{ | c | c | c | c | c | }
\hline
Ref. & $\left(\bm{A}_{t},\bm{B}_{t}\right)$ & $\bm{Q}_{t}\succeq\bm{0}$ & Distributions $\mu_0,\mu_1$ & Markov kernel $\kappa$\\[0.5ex] 
 \hline\hline
\cite{Sch31,Sch32} & $\left(\bm{0},\bm{I}\right)$ & $\bm{0}$ & generic & \eqref{HeatKernel}\\[0.5ex]
\hline
\cite{chenPartI} & generic & $\bm{0}$ & Gaussian & \eqref{LinearKernel}\\[0.5ex]
\hline
\cite{chen2015optimal,chenPartIII} & generic & generic & Gaussian & n/a\\[0.5ex]
\hline
\cite{teter2023contraction} & generic & $\bm{0}$ & generic & \eqref{LinearKernel}\\[0.5ex]
\hline
\cite{teter2024schr,teter2024weyl} & $\left(\bm{0},\bm{I}\right)$ & fixed $\bm{Q}\succeq\bm{0}$ & generic & \cite[(A.22)]{teter2024schr}, \cite[(43)]{teter2024weyl}\\[0.5ex]
\hline
This work & generic & generic & generic & \eqref{KernelInTermsOfM}\\[0.5ex]
\hline
\end{tabular}
\caption{Comparison of related works on LQ SBP.}
\label{table:relatedworksLQSBP}
\end{table}

\subsubsection*{Organization}
In Sec. \ref{sec:kernelfromdist}, we motivate the postulated form (see \eqref{kappaExpDist}) for the Markov kernel of interest in terms of a distance function. In Sec. \ref{sec:distfromocp}, we re-visit the known Markov kernels to demonstrate that the distance function appearing in our postulated form can be obtained by solving an associated deterministic optimal control problem. 

Motivated by the structural observations made in Sections \ref{sec:kernelfromdist} and \ref{sec:distfromocp}, we next follow the computational template: Markov kernel $\longleftarrow$ distance function $\longleftarrow$ deterministic optimal control problem, to derive the Markov kernel for the It\^{o} diffusion \eqref{LinearSDE} with rate of killing $\frac{1}{2}\bm{x}^{\top}\bm{Q}_{t}\bm{x}$ for $\bm{Q}_{t}\succeq\bm{0}$. In particular, the corresponding distance is derived in Sec. \ref{subsec:gettingdist}, which we then use to compute the Markov kernel of interest in Sec. \ref{subsec:gettingkappa}. Our main result is Theorem \ref{thm:kernelLQ}. In Sec. \ref{subsec:FromkappatoBridge}, we explain how the derived kernel helps to solve the LQ SBP with non-Gaussian endpoints. Sec. \ref{sec:conclusions} concludes the work.


\section{Markov Kernels and Distances}\label{sec:kernelfromdist}
Both \eqref{HeatKernel} and \eqref{LinearKernel} are of the form
\begin{align}
\kappa = c\left(t,t_{0}\right)\exp\left(-\dfrac{1}{2}\dist_{tt_{0}}^{2}\left(\bm{x},\bm{y}\right)\right)
\label{kappaExpDist}    
\end{align}
for some suitable distance function $$\dist_{tt_{0}}:\mathbb{R}^{n}\times\mathbb{R}^{n}\mapsto\mathbb{R}_{\geq 0}.$$
The subscript $tt_{0}$ signifies the distance function's parametric dependence on $t,t_{0}$.

Because \eqref{HeatKernel} and \eqref{LinearKernel} are both transition probabilities, the distance function $\dist_{tt_{0}}$ uniquely determines $\kappa$. Once $\dist_{tt_{0}}$ is identified, the pre-factor $c\left(t,t_{0}\right)$, and hence $\kappa$ is completely determined by the normalization condition $\int_{\mathbb{R}^{n}}\kappa\differential\bm{y} = 1$.

Interestingly, the $\kappa$ associated with problem \eqref{ClassicalSBPstatecost} derived in \cite{teter2024schr,teter2024weyl} is also of the form \eqref{kappaExpDist} even though $\kappa$ then is not a transition probability, and $c(t,t_0)$ does not follow from normalization. 

These exemplars hint that for controlled dynamics over $\mathbb{R}^{n}$, the distance function might be induced by some principle of least action.

That Markov kernels are related to distances is in itself not a new observation. The most well-known result relating kernels and distances is \emph{Varadhan's formula} \cite{varadhan1967behavior,molchanov1975diffusion,crane2017heat} which says that the heat kernel $\kappa_{\mathrm{Heat}}^{\mathcal{M}}$ on a complete Riemannian manifold $\mathcal{M}$ satisfies
\begin{align}
\displaystyle\lim_{t\downarrow t_{0}}t\log \kappa_{\mathrm{Heat}}^{\mathcal{M}}\left(t_{0},\bm{x},t,\bm{y}\right) = -\frac{1}{2}\dist^{2}\left(\bm{x},\bm{y}\right)
\label{VaradhanFormula}	
\end{align}
uniformly on every compact subsets of $\mathcal{M}\times\mathcal{M}$, where $\dist$ is the minimal geodesic distance connecting $\bm{x},\bm{y}\in\mathcal{M}$. In other words, $\dist$ can be recovered as the short time asymptotic of the heat kernel on $\mathcal{M}$. In the context of Varadhan's formula \eqref{VaradhanFormula}, the $L$ in \eqref{KolForwardPDEIVP} is the Laplace-Beltrami operator on $\mathcal{M}$. 

For specific manifolds, few exact or asymptotic relations are also known \cite[Ch. VI]{chavel1984eigenvalues}, \cite{grigor1998heat}, \cite[Ch. 5]{hsu2002stochastic}:
\begin{subequations}
\begin{align}
\kappa_{\mathrm{Heat}}^{\mathbb{H}^{3}}\left(t_{0},\bm{x},t,\bm{y}\right) = &\dfrac{1}{\left(4\pi(t-t_{0})\right)^{3/2}}\dfrac{\dist(\bm{x},\bm{y})}{\sinh\dist(\bm{x},\bm{y})}\nonumber\\
&\times \exp\left(-t-\dfrac{\dist^{2}(\bm{x},\bm{y})}{4(t-t_0)}\right), \label{HeatH3Exact}\\
\kappa_{\mathrm{Heat}}^{\mathbb{H}^{n}}\left(t_{0},\bm{x},t,\bm{y}\right) \sim &\dfrac{1}{\left(4\pi(t-t_{0})\right)^{n/2}}\!\left(\!\dfrac{\dist(\bm{x},\bm{y})}{\sinh\dist(\bm{x},\bm{y})}\!\right)^{\!(n-1)/2}\nonumber\\
&\times \exp\left(-\dfrac{\dist^{2}(\bm{x},\bm{y})}{4(t-t_0)}\right), \label{HeatHnAsymp}\\
\kappa_{\mathrm{Heat}}^{\mathbb{S}^{n}}\left(t_{0},\bm{x},t,\bm{y}\right) \sim &\dfrac{1}{\left(4\pi(t-t_{0})\right)^{n/2}}\!\left(\!\dfrac{\dist(\bm{x},\bm{y})}{\sin\dist(\bm{x},\bm{y})}\!\right)^{\!(n-1)/2}\nonumber\\
&\times \exp\left(-\dfrac{\dist^{2}(\bm{x},\bm{y})}{4(t-t_0)}\right), \label{HeatSnAsymp}
\end{align}
\label{HeatKernelExactOrAsymp}
\end{subequations}
where $\mathbb{H}^{n}, \mathbb{S}^{n}$ denote the $n$ dimensional hyperbolic manifold and the sphere, respectively. The symbol $\sim$ in \eqref{HeatKernelExactOrAsymp} denotes asymptotic equivalence.

In contrast, the form \eqref{kappaExpDist} that we focus on is rather specific, and is for a flat geometry. It is then natural to speculate that $\dist_{tt_{0}}$ may arise from the finite horizon reachability constraint over $[t_{0},t]$ imposed by the controlled dynamics, i.e., $\dist_{tt_{0}}$ is of \emph{sub-Riemannian} or \emph{Carnot-Carath{\'e}odory} type \cite{brockett1982control,strichartz1986sub,gromov1996carnot}. This motivates us to formulate $\dist_{tt_{0}}$ as the minimal value of an action integral explained next.


\section{Distances from Optimal Control}\label{sec:distfromocp}
We postulate that for $\kappa$ of the form \eqref{kappaExpDist}, the $\dist_{tt_{0}}$ is induced by a \emph{deterministic optimal control problem (OCP)} of particular structure. The objective for this OCP is
\begin{align}
\dfrac{1}{2}\dist_{tt_{0}}^{2}(\bm{x},\bm{y}) = \underset{\bm{u}_{\tau}}{\min}\int_{t_{0}}^{t}\left(\frac{1}{2}\vert\bm{u}_{\tau}\vert^{2} + q(\bm{z}_{\tau})\right)\differential\tau.
\label{detOCP}    
\end{align}
The constraint for this OCP is a \emph{controlled ODE} obtained by replacing $\differential\bm{w}_{\tau}$ in the underlying It\^{o} diffusion with a controlled drift $\bm{u}_{\tau}\:\differential \tau$, and boundary conditions $\bm{z}_{\tau}(\tau=t_0)=\bm{x}$, $\bm{z}_{\tau}(\tau=t) = \bm{y}$. Notice that in objective \eqref{detOCP}, the state cost $q$ is the rate of creation or killing of probability mass. Let us verify this postulate for known cases discussed before.\\

\noindent\textbf{Heat kernel \eqref{HeatKernel}.} Here $q\equiv 0$, and $\differential\bm{w}_{\tau}\mapsto \bm{u}\:\differential \tau$ in \eqref{BrownianSDE} gives the controlled ODE $\dot{\bm{z}}_{\tau} = \sqrt{2}\:\bm{u}_{\tau}$ where dot denotes derivative w.r.t. $\tau\in[t_0,t]$. This leads to the deterministic OCP
\begin{subequations}
\begin{align}
\dfrac{1}{2}\dist_{tt_{0}}^{2}(\bm{x},\bm{y}) = &\underset{\bm{u}_{\tau}}{\min}\int_{t_{0}}^{t}\frac{1}{2}\vert\bm{u}_{\tau}\vert^{2} \differential\tau\label{detOCPHeatObj}\\
& \dot{\bm{z}}_{\tau} = \sqrt{2}\:\bm{u}_{\tau},\label{detOCPHeatODE}\\ & \bm{z}_{\tau}(\tau=t_0)=\bm{x}, \quad \bm{z}_{\tau}(\tau=t) = \bm{y}.\label{detOCPHeatBC}
\end{align}
\label{detOCPHeat}
\end{subequations}
For solving \eqref{detOCPHeat}, we apply the Pontryagin's minimum principle to get the optimal control $\bm{u}_{\tau}^{\mathrm{opt}} = -\sqrt{2}\bm{\lambda}_{\tau}^{\mathrm{opt}}$ where $\bm{\lambda}_{\tau}^{\mathrm{opt}}$ is the optimal costate. The optimal state $\ddot{\bm{z}}_{\tau}^{\mathrm{opt}} = 0$, or $\dot{\bm{z}}_{\tau}^{\mathrm{opt}} = -2\bm{\lambda}_{\tau}^{\mathrm{opt}} = \bm{\alpha}$, or equivalently $\bm{z}_{\tau}^{\mathrm{opt}} = \bm{\alpha} \tau + \bm{\beta}$ for some constant $\bm{\alpha},\bm{\beta}\in\mathbb{R}^{n}$. 

Using \eqref{detOCPHeatBC}, we find $\bm{\alpha} = (\bm{x}-\bm{y})/(t_0 - t)$, and $\bm{\lambda}_{\tau}^{\mathrm{opt}} = -\bm{\alpha}/2 = (\bm{x}-\bm{y})/2(t-t_0)$. Hence the optimal value in \eqref{detOCPHeatObj}:
$$\dfrac{1}{2}\dist_{tt_{0}}^{2}(\bm{x},\bm{y}) = \frac{1}{2}\cdot 2\vert \bm{\lambda}_{\tau}^{\mathrm{opt}}\vert^2 \cdot (t-t_{0}) = \dfrac{\vert \bm{x} - \bm{y}\vert^2}{4(t-t_{0})},$$
which indeed transcribes \eqref{HeatKernel} into the form \eqref{kappaExpDist}. The normalization condition for $\kappa$ leads to evaluating a Gaussian integral, which determines the pre-factor $$c(t,t_{0}) = (4\pi(t-t_0))^{-n/2}.$$\\

\noindent\textbf{Linear kernel \eqref{LinearKernel}.} Here $q\equiv 0$, and $\differential\bm{w}_{\tau}\mapsto \bm{u}\:\differential \tau$ in \eqref{LinearSDE} yields the deterministic OCP
\begin{subequations}
\begin{align}
\dfrac{1}{2}\dist_{tt_{0}}^{2}(\bm{x},\bm{y}) = &\underset{\bm{u}_{\tau}}{\min}\int_{t_{0}}^{t}\frac{1}{2}\vert\bm{u}_{\tau}\vert^{2} \differential\tau\label{detOCPLinearObj}\\
& \dot{\bm{z}}_{\tau} = \bm{A}_{\tau}\bm{z}_{\tau} + \sqrt{2}\:\bm{B}_{\tau}\bm{u}_{\tau},\label{detOCPLinearODE}\\ & \bm{z}_{\tau}(\tau=t_0)=\bm{x}, \quad \bm{z}_{\tau}(\tau=t) = \bm{y}.\label{detOCPLinearBC}
\end{align}
\label{detOCPLinear}
\end{subequations}
Problem \eqref{detOCPLinear} is that of minimum effort state steering for controllable linear system, and its solution is commonplace in optimal control textbooks; see e.g., \cite[p. 194]{lee1967foundations}. The optimal value in \eqref{detOCPLinearObj}:
$$\dfrac{1}{2}\dist_{tt_{0}}^{2}(\bm{x},\bm{y}) = \frac{\left(\boldsymbol{\Phi}_{tt_0} \boldsymbol{x}-\boldsymbol{y}\right)^{\top} \boldsymbol{\Gamma}_{t t_0}^{-1}\left(\boldsymbol{\Phi}_{tt_0} \boldsymbol{x}-\boldsymbol{y}\right)}{4\left(t-t_0\right)},$$
which indeed transcribes \eqref{LinearKernel} into the form \eqref{kappaExpDist}. As before, the normalization condition for $\kappa$ determines the pre-factor $$c(t,t_{0}) = \left(4 \pi\left(t-t_0\right)\right)^{-n/2} \operatorname{det}\left(\boldsymbol{\Gamma}_{tt_0}\right)^{-1 / 2}$$ via Gaussian integral.\\

\noindent\textbf{The kernel in \cite{teter2024schr,teter2024weyl}.} The kernel derived in \cite{teter2024schr,teter2024weyl} corresponds to the It\^{o} diffusion \eqref{BrownianSDE} with reaction rate $\frac{1}{2}\bm{z}^{\top}\bm{Q}\bm{z}$, $\bm{Q}\succeq\bm{0}$. Without loss of generality \cite[Sec. 4.1]{teter2024schr}, we consider state coordinates where a given $\bm{Q}\succ\bm{0}$ is diagonalized to yield the diagonal matrix $2\bm{D}\succ\bm{0}$. The factor $2$ scaling is unimportant but kept for consistency with \cite[Sec. 4]{teter2024schr}. Importantly, the Laplacian operator is invariant under this change of coordinates. For the case $\bm{Q}\succeq\bm{0}$, we refer the readers to \cite[Sec. 4.3]{teter2024schr}.

In the new state coordinates, $q\equiv \frac{1}{2}\bm{z}_{\tau}^{\top}2\bm{D}\bm{z}_{\tau}$, and $\differential\bm{w}_{\tau}\mapsto \bm{u}\:\differential \tau$ in \eqref{BrownianSDE} gives the following modified version of the OCP \eqref{detOCPHeat}:
\begin{subequations}
\begin{align}
\dfrac{1}{2}\dist_{tt_{0}}^{2}(\bm{x},\bm{y}) = &\underset{\bm{u}_{\tau}}{\min}\!\int_{t_{0}}^{t}\!\!\left(\frac{1}{2}\vert\bm{u}_{\tau}\vert^{2} + \frac{1}{2}\bm{z}_{\tau}^{\top}2\bm{D}\bm{z}_{\tau}\!\right) \differential\tau\label{detOCPsiconObj}\\
& \dot{\bm{z}}_{\tau} = \sqrt{2}\:\bm{u}_{\tau},\label{detOCPsiconODE}\\ & \bm{z}_{\tau}(\tau=t_0)=\bm{x}, \quad \bm{z}_{\tau}(\tau=t) = \bm{y}.\label{detOCPsiconBC}
\end{align}
\label{detOCPsicon}
\end{subequations}
Applying Pontryagin's minimum principle, we get the optimal control $\bm{u}_{\tau}^{\mathrm{opt}} = -\sqrt{2}\bm{\lambda}_{\tau}^{\mathrm{opt}}$, the optimal costate ODE $\dot{\bm{\lambda}}_{\tau}^{\mathrm{opt}} = -2\bm{D}\bm{z}_{\tau}^{\mathrm{opt}}$, and the optimal state ODE 
\begin{align}
\ddot{\bm{z}}_{\tau}^{\mathrm{opt}} = -2\dot{\bm{\lambda}}_{\tau}^{\mathrm{opt}} = 4\bm{D}\bm{z}_{\tau}^{\mathrm{opt}}.
\label{ddotzSICON}    
\end{align}
Letting $\omega_i := 2\sqrt{{D}_{ii}}$ for all $i=1,\hdots,n$, we solve for the components of \eqref{ddotzSICON} as 
\begin{align}
\left(z_{\tau}^{\mathrm{opt}}\right)_{i}=a_{i}e^{\omega_i\tau}+b_{i}e^{-\omega_i\tau},
\label{zicomponentsol}    
\end{align}
and use \eqref{detOCPsiconBC} to determine the constants
\begin{align}
a_i = \dfrac{-x_i e^{-\omega_i t}+y_i e^{-\omega_i t_0}}{2 \sinh \left(\omega_i\left(t-t_0\right)\right)},
b_i = \dfrac{x_i e^{\omega_i t}-y_i e^{\omega_i t_0}}{2 \sinh \left(\omega_i\left(t-t_0\right)\right)}. 
\label{aibi}
\end{align}
Hence the optimal value in \eqref{detOCPsiconObj}:
\begin{align}
&\dfrac{1}{2}\dist_{tt_{0}}^{2}(\bm{x},\bm{y})\nonumber\\
=& \!\!\displaystyle\int_{t_{0}}^{t}\!\!\left(\dfrac{1}{2}\vert\bm{u}_{\tau}^{\mathrm{opt}}\vert^2 + \dfrac{1}{4}\displaystyle\sum_{i=1}^{n}\omega_{i}^{2}\left(\left(z_{\tau}^{\mathrm{opt}}\right)_{i}\right)^{2}\right)\differential\tau\nonumber\\
=& \!\!\displaystyle\int_{t_{0}}^{t}\!\!\left(\dfrac{1}{4}\vert\bm{z}_{\tau}^{\mathrm{opt}}\vert^2 + \dfrac{1}{4}\displaystyle\sum_{i=1}^{n}\omega_{i}^{2}\left(\left(z_{\tau}^{\mathrm{opt}}\right)_{i}\right)^{2}\right)\differential\tau\nonumber\\
=&\dfrac{1}{2}\displaystyle\sum_{i=1}^{n}\dfrac{\omega_{i}\left(x_i^2 + y_i^2\right)\cosh\left(\omega_{i}(t-t_{0})\right) - 2\omega_i x_i y_i}{2\sinh\left(\omega_{i}(t-t_{0})\right)},\label{halfdistsquaredSICONkernel} 
\end{align}
where we have used \eqref{zicomponentsol}-\eqref{aibi} followed by algebraic simplification using the identity $\sinh(2(\cdot))=2\sinh(\cdot)\cosh(\cdot)$.

We note that the expression \eqref{halfdistsquaredSICONkernel} indeed transcribes the Markov kernel in \cite[Eq. (43)]{teter2024weyl} or that in \cite[Eq. (A.22)]{teter2024schr} into the form \eqref{kappaExpDist}. Since this Markov kernel is not a transition probability, the normalization condition no longer holds, and the pre-factor
\begin{align}
c(t,t_0)=\displaystyle\prod_{i=1}^{n}\left(\dfrac{\omega_i}{4\pi\sinh\left(\omega_{i}(t-t_{0})\right)}\right)^{1/2}
\label{PrefactorKernelSICON}   
\end{align}
does not follow from there. However, having determined \eqref{halfdistsquaredSICONkernel}, the pre-factor \eqref{PrefactorKernelSICON} can be obtained by substituting \eqref{kappaExpDist} with \eqref{halfdistsquaredSICONkernel} in \eqref{ReactionAdvectionDiffusionPDEIVP}. See Appendix \ref{DerivationOfcClassicalReactionDiffusion} for details of this simple but non-trivial computation.


\section{Distance and kernel for LQ Non-Gaussian Distribution Steering}\label{sec:MainSection}
Having seen three exemplars for the computational pipeline:
$$\text{Markov kernel}\;\kappa\longleftarrow\dist\longleftarrow\text{deterministic OCP},$$
we now apply the same to derive $\kappa$ for the It\^{o} diffusion \eqref{LinearSDE} with reaction rate $q(\bm{z})\equiv\frac{1}{2}\bm{z}^{\top}\bm{Q}_{\tau}\bm{z}$, $\tau\in[t_0,t]$. For convenience, let $\hat{\bm{B}}_{\tau}:= \sqrt{2}\:\bm{B}_{\tau}$. 

Our standing assumptions are:\\
\noindent\textbf{A1.} $\left(\bm{A}_\tau,\bm{B}_\tau\right)$, and thus $\left(\bm{A}_\tau,\hat{\bm{B}}_\tau\right)$, is controllable matrix-valued trajectory pair in $\left(\mathbb{R}^{n\times n},\mathbb{R}^{n\times m}\right)$ that is bounded and continuous in $\tau\in[t_0,t]$.

\noindent\textbf{A2.} the matricial trajectory $\bm{Q}_{\tau}\succeq\bm{0}$ is continuous and bounded w.r.t. $\tau\in[t_0,t]$, and $\bm{Q}_{s}\succ\bm{0}$ for some $s\in[t_0,t]$.

Finding the Markov kernel $\kappa$ in this setting enables the solution of the LQ SBP with generic \emph{non-Gaussian endpoint distributions} having finite second moments, i.e., the solution of the problem:
\begin{subequations}
\begin{align}
&\underset{\left(\mu^{\bm{u}},\bm{u}\right)}{\inf}\int_{\mathbb{R}^{n}}\int_{t_{0}}^{t_{1}}\bigg\{\frac{1}{2}\vert\bm{u}\vert^{2} + \frac{1}{2}\left(\bm{x}_t^{\bm{u}}\right)^{\top}\bm{Q}_{t}\bm{x}_t^{\bm{u}}\bigg\}\differential t \:\differential\mu^{\bm{u}}(\bm{x}_{t}^{\bm{u}})\label{LinearSBPstatecostObj}\\
&\text{subject to}\quad \differential\bm{x}_{t}^{\bm{u}} = \left(\bm{A}_t\bm{x}_t^{\bm{u}} +\bm{B}_t\bm{u}_{t}\right)\:\differential t + \sqrt{2}\bm{B}_{t}\:\differential\bm{w}_{t},\label{LinearSBPstatecostSDE}\\
&\qquad\qquad\quad\bm{x}_{t}^{\bm{u}}(t=t_0) \sim \mu_{0}, \quad \bm{x}_{t}^{\bm{u}}(t=t_1) \sim \mu_{1}. \label{LinearSBPstatecostEndpointConstr}\end{align}
\label{LinearSBPstatecost}
\end{subequations}

For $\mu_0,\mu_1$ with finite second moments, the existence-uniqueness for the solution of \eqref{LinearSBPstatecost} is guaranteed. This follows from transcribing \eqref{LinearSBPstatecost} to a stochastic calculus of variations problem involving the relative entropy a.k.a. the Kullback-Leibler divergence $\mathrm{D}_{\mathrm{KL}}\left(\cdot\parallel\cdot\right)$ minimization:
\begin{align}
\underset{\mathbb{P}\in\Pi_{01}}{\min}\quad{\mathrm{D}_{\mathrm{KL}}}\left(\mathbb{P}\parallel \dfrac{\exp\left(-\frac{1}{4}\int_{t_0}^{t_1}\bm{x}^{\top}\bm{Q}_{t}\bm{x}\:\differential t\right)\mathbb{W}}{Z}\right)
\label{LDPforLQSBP}    
\end{align}
where
\begin{align}
&\Pi_{01} := \{\mathbb{M}\in\mathcal{M}\left(\mathcal{C}\left(\left[t_0,t_{1}\right];\mathbb{R}^{n}\right)\right)\mid\mathbb{M}\;\text{has marginal}\;\mu_{k}\nonumber\\
&\qquad\qquad\qquad\qquad\qquad\qquad\text{at time}\;t_{k}\:\forall k\in\{0,1\}\},
\label{defPi01}    
\end{align}
the $\mathcal{M}\left(\mathcal{C}\left(\left[t_0,t_{1}\right];\mathbb{R}^{n}\right)\right)$ is the set of probability measures on the path space $\mathcal{C}\left(\left[t_0,t_{1}\right];\mathbb{R}^{n}\right)$ generated by the It\^{o} diffusion \eqref{LinearSDE}, the $\mathbb{W}\in\mathcal{M}\left(\mathcal{C}\left(\left[t_0,t_{1}\right];\mathbb{R}^{n}\right)\right)$ is the Wiener measure, and the $Z$ is a normalizing constant. Problem \eqref{LDPforLQSBP} seeks to find the most likely measure-valued path generated by \eqref{LinearSDE} w.r.t. a weighted Wiener (i.e., Gibbs) measure defined by the quadratic state cost, with endpoint marginal constraints given by the problem data $\mu_0,\mu_1$. For the equivalence of the formulations \eqref{LinearSBPstatecost} and \eqref{LDPforLQSBP}, see e.g., \cite{wakolbinger1990schrodinger,dawson1990schrodinger,leonard2011stochastic}, \cite[Sec. 4.1]{teter2025probabilistic}. The existence-uniqueness of solution for \eqref{LDPforLQSBP} is a consequence of strict convexity of $\mathrm{D}_{\mathrm{KL}}\left(\cdot\parallel\cdot\right)$ w.r.t. the first argument.

Unlike the three kernels in Sec. \ref{sec:distfromocp}, the $\kappa$ for problem \eqref{LinearSBPstatecost} is not known. Moreover, previously mentioned approaches (Hermite polynomials, Weyl calculus) to compute $\kappa$ that works for simpler cases no longer generalize in any obvious way. We postulate that the form \eqref{kappaExpDist} for the Markov kernel holds here as well. In Sections \ref{subsec:gettingdist}-\ref{subsec:gettingkappa}, we will verify the same. In Section \ref{subsec:FromkappatoBridge}, we will explain how the resulting kernel helps to solve the LQ SBP \eqref{LinearSBPstatecost}.

\subsection{Getting $\dist$}\label{subsec:gettingdist}
Inspired by the instances in Sec. \ref{sec:distfromocp}, we define $\dist$ through the following deterministic OCP:
\begin{subequations}
\begin{align}
\dfrac{1}{2}\dist_{tt_{0}}^{2}(\bm{x},\bm{y}) = &\underset{\bm{u}_{\tau}}{\min}\int_{t_{0}}^{t}\left(\frac{1}{2}\vert\bm{u}_{\tau}\vert^{2}+\frac{1}{2}\bm{z}_{\tau}^{\top}\bm{Q}_{\tau}\bm{z}_{\tau}\right) \differential\tau\label{detOCPLinearQuadraticObj}\\
& \dot{\bm{z}}_{\tau} = \bm{A}_{\tau}\bm{z}_{\tau} + \sqrt{2}\:\bm{B}_{\tau}\bm{u}_{\tau},\label{detOCPLinearQuadraticODE}\\ & \bm{z}_{\tau}(\tau=t_0)=\bm{x}, \quad \bm{z}_{\tau}(\tau=t) = \bm{y}.\label{detOCPLinearQuadraticBC}
\end{align}
\label{detOCPLinearQuadratic}
\end{subequations}
This is an atypical linear quadratic OCP in that when $\bm{x},\bm{y}$ are not too close to the origin, the optimal state trajectory trades off the soft penalty (state cost-to-go) in deviating away from origin with that of meeting the hard endpoint constraints \eqref{detOCPLinearQuadraticBC} within the given deadline.

We need the following result from \cite[p. 140-141]{brockett1970finite} paraphrased below to suit our notations.
\begin{proposition}\cite[p. 140-141]{brockett1970finite}\label{Prop:etaetahat}
Suppose there exists symmetric matrix $\bm{K}_1$ such that the solution map\footnote{The mapping $\bm{\Pi}(\tau, \bm{K}_1, t)$ is understood as the solution of \eqref{RiccatiODE} at any $\tau\in[t_0,t]$ solved backward in time with initial condition \eqref{RiccatiIC}.} $\bm{\Pi}(\tau, \bm{K}_1, t)$ of the Riccati matrix ODE initial value problem
\begin{subequations}
\begin{align}
    &\dot{\bm{K}}_{\tau} = - \bm{A}^{\top}_{\tau} \bm{K}_{\tau} - \bm{K}_{\tau}\bm{A}_{\tau} + \bm{K}_{\tau}\hat{\bm{B}}_{\tau}\hat{\bm{B}}^{\top}_{\tau}\bm{K}_\tau - \bm{Q}_\tau,\label{RiccatiODE}\\
    &\bm{K}_{\tau = t} = \bm{K}_1,\label{RiccatiIC}
\end{align}
\label{RiccatiIVP}
\end{subequations}
exists for all $\tau\in[t_0, t]$. Let 
\begin{align}
\hat{\bm{A}}_{\tau} := \bm{A}_{\tau} - \hat{\bm{B}}_{\tau}\hat{\bm{B}}^{\top}_{\tau} \bm{\Pi}(\tau, \bm{K}_1, t) \quad \forall\tau\in[t_0,t].
\label{defAhat}    
\end{align}
Then the differentiable trajectory $\bm{z}_{\tau}$ minimizes
\begin{align*}
    \eta =  &\int_{t_{0}}^{t}\left(\frac{1}{2}\vert\bm{u}_{\tau}\vert^{2}+\frac{1}{2}\bm{z}_{\tau}^{\top}\bm{Q}_{\tau}\bm{z}_{\tau}\right) \differential\tau \\
    \text{subject to} \quad &\dot{\bm{z}}_{\tau} = \bm{A}_{\tau}\bm{z}_{\tau} + \hat{\bm{B}}_{\tau}\bm{u}_{\tau},\\
    &\bm{z}_{\tau = t_0} = \bm{x}, \quad \bm{z}_{\tau = t} = \bm{y},
\end{align*}
if and only if it also minimizes
\begin{align*}
    \hat{\eta} =  &\int_{t_{0}}^{t}\frac{1}{2}\vert\bm{v}_{\tau}\vert^{2} \differential\tau \\
    \text{subject to} \quad &\dot{\bm{z}}_{\tau} = \hat{\bm{A}}_{\tau}\bm{z}_{\tau} + \hat{\bm{B}}_{\tau}\bm{v}_{\tau},\\
    &\bm{z}_{\tau = t_0} = \bm{x}, \quad \bm{z}_{\tau = t} = \bm{y}.
\end{align*}
Additionally, along any trajectory satisfying the boundary conditions, we have
\begin{align}
    \eta = \hat{\eta} + \frac{1}{2}\bm{x}^{\top}\bm{\Pi}(t_0, \bm{K}_1, t)\bm{x} - \frac{1}{2}\bm{y}^{\top}\bm{K}_1 \bm{y}.
    \label{CostRelation}
\end{align}
\end{proposition}
\begin{remark}[Existence, uniqueness and positive semi-definiteness of $\bm{\Pi}$]\label{Remark:Controllability}
Proposition \ref{Prop:etaetahat} does not assume the controllability of \eqref{detOCPLinearQuadraticODE}. However, under our standing controllability assumption \textbf{A1} and the positive semi-definiteness assumption for $\bm{Q}_{\tau}$ $\forall\tau\in[t_0,t]$ in \textbf{A2}, the existence-uniqueness of the mapping $\bm{\Pi}$ is guaranteed for all $\tau\in[t_0,t]$. Furthermore, the unique solution $\bm{\Pi}(\tau, \bm{K}_1, t)\succeq\bm{0}$ for any $\bm{K}_{1}\succeq\bm{0}$. See e.g., \cite[Thm. 8-10]{kuvcera1973review}.
\end{remark}
\begin{remark}[Strict positive definiteness of $\bm{\Pi}$]\label{Remark:PosDefPi}
The assumption $\bm{Q}_{s}\succ\bm{0}$ for some $s\in[t_0,t]$, stated in the latter part of $\textbf{A2}$, further ensures that $\bm{\Pi}(\tau, \bm{K}_1, t)\succ\bm{0}$ for any $\bm{K}_{1}\succeq\bm{0}$. This is a consequence of the matrix variations of constants formula \cite[p. 59, exercise 1 in p. 162]{brockett1970finite}, \cite[Sec. II]{porter1967matrix}.
\end{remark}


Since the optimal $\hat{\eta}$ in Proposition \ref{Prop:etaetahat} is the cost for minimum energy state transfer from $\bm{x}$ at $t_0$ to $\bm{y}$ at $t$ over the LTV system matrix pair $(\hat{\bm{A}}_{\tau},\hat{\bm{B}}_{\tau})$, we have
\begin{align}
\hat{\eta}_{\mathrm{opt}} = \dfrac{1}{2}\left(\hat{\bm{\Phi}}_{t t_0}\bm{x} - \bm{y}\right)^{\top} \hat{\bm{\Gamma}}_{t t_0}^{-1}\left(\hat{\bm{\Phi}}_{t t_0}\bm{x} - \bm{y}\right),
\label{etahatopt}    
\end{align}
where $\hat{\bm{\Phi}}_{t t_0}$ is the state transition matrix from $t_0$ to $t$ for the state matrix $\hat{\bm{A}}_\tau$ in \eqref{defAhat}. Likewise, $\hat{\bm{\Gamma}}_{t t_0}$ in \eqref{etahatopt} is the controllability Gramian for the pair $(\hat{\bm{A}}_{\tau},\hat{\bm{B}}_{\tau})$. 

To see why \eqref{etahatopt} is well-defined, recall that for linear systems, controllability remains invariant under state feedback. Specifically, the following Proposition \ref{prop:controllability} holds. For its proof in the LTV setting, see \cite{viswanadham1975invariance}; the proof in the linear time-invariant setting appeared earlier in \cite[Thm. 3]{brockett1965poles} .
\begin{proposition}\cite{viswanadham1975invariance}\label{prop:controllability}
Let $\hat{\bm{A}}_\tau$ be given by \eqref{defAhat}. If the LTV system defined by the pair $(\bm{A}_{\tau},\hat{\bm{B}}_{\tau})$ is controllable, then so is the LTV system defined by the pair $(\hat{\bm{A}}_{\tau},\hat{\bm{B}}_{\tau})$.     
\end{proposition}
From Proposition \ref{prop:controllability}, it follows that $\bm{\Gamma}_{t t_0}\succ\bm{0}$ implies $\hat{\bm{\Gamma}}_{t t_0}\succ\bm{0}$. In particular, $\hat{\bm{\Gamma}}_{t t_0}$ is non-singular, and \eqref{etahatopt} is well-defined. 

Since the optimal $\eta$ in \eqref{CostRelation}, and thus the optimal value in \eqref{detOCPLinearQuadraticObj} must be independent of $\bm{K}_{1}$, we can set $\bm{K}_{1}\equiv\bm{0}$ without loss of generality. Combining this observation with \eqref{etahatopt}, we find the optimal value in \eqref{detOCPLinearQuadraticObj} as
\begin{align}
\dfrac{1}{2}\dist_{tt_{0}}^{2}(\bm{x},\bm{y})= \!\dfrac{1}{2}\!\begin{pmatrix}
\bm{x}\\ \bm{y}
\end{pmatrix}^{\!\!\top}\bm{M}_{tt_{0}}\begin{pmatrix}
\bm{x}\\ \bm{y}
\end{pmatrix},
\label{HalfDistSquaredABQ}
\end{align}
where
\begin{align}
\bm{M}_{tt_0}\! := \!\begin{bmatrix}
\hat{\bm{\Phi}}_{tt_0}^{\top}\hat{\bm{\Gamma}}_{tt_0}^{-1}\hat{\bm{\Phi}}_{tt_0}+\bm{\Pi}(t_0, \bm{0}, t) & -\hat{\bm{\Phi}}_{tt_0}^{\top}\hat{\bm{\Gamma}}_{tt_0}^{-1}\\
&\\
-\hat{\bm{\Gamma}}_{tt_0}^{-1}\hat{\bm{\Phi}}_{tt_0} & \hat{\bm{\Gamma}}_{tt_0}^{-1}
\!\end{bmatrix}\!.
\label{defMtt0}
\end{align}

Using the Schur complement lemma, we can check that $\bm{M}_{tt_0}\succ\bm{0}$, as expected. Specifically, from \eqref{defMtt0}, note that
\begin{subequations}
\begin{align}
&\hat{\bm{\Gamma}}_{tt_0}^{-1}\succ\bm{0},\\
& \hat{\bm{\Phi}}_{tt_0}^{\top}\hat{\bm{\Gamma}}_{tt_0}^{-1}\hat{\bm{\Phi}}_{tt_0}+\bm{\Pi}(t_0, \bm{0}, t) - \left(-\hat{\bm{\Phi}}_{tt_0}^{\top}\hat{\bm{\Gamma}}_{tt_0}^{-1}\right)\hat{\bm{\Gamma}}_{tt_0}\left(-\hat{\bm{\Gamma}}_{tt_0}^{-1}\hat{\bm{\Phi}}_{tt_0}\right)\nonumber\\
&= \bm{\Pi}(t_0, \bm{0}, t)\succ\bm{0},
\end{align}
\label{SchurComplementOfMwrtRightBottomBlock}    
\end{subequations}
thanks to Remark \ref{Remark:PosDefPi}.

The formulae \eqref{HalfDistSquaredABQ}-\eqref{defMtt0} generalize the previously known result in \cite{teter2024schr,teter2024weyl} for zero prior dynamics. In particular, the previously known result can be recovered from \eqref{HalfDistSquaredABQ}-\eqref{defMtt0} as follows (proof in Appendix \ref{AppProofOfProp:SpecialCaseAzeroBidentity}). 
\begin{proposition}\label{prop:SpecialCaseAzeroBidentity}
When $\bm{Q}_{\tau}=2\bm{D}\succ\bm{0}$ (constant positive diagonal matrix) and $(\bm{A}_{\tau},\hat{\bm{B}}_{\tau})\equiv(\bm{0},\sqrt{2}\bm{I})$ $\forall \tau\in[t_0,t]$, the formulas \eqref{HalfDistSquaredABQ}-\eqref{defMtt0} reduce to \eqref{halfdistsquaredSICONkernel}.    
\end{proposition}



\subsection{Getting $\kappa$}\label{subsec:gettingkappa}
Now that we have determined the distance functional in the generic LQ setting given by \eqref{HalfDistSquaredABQ}-\eqref{defMtt0}, what remains in finding $\kappa$ in the postulated form \eqref{kappaExpDist}, i.e., in the form
\begin{align}
\kappa\left(t_0,\bm{x},t,\bm{y}\right) = c(t,t_0)\exp\left(- \frac{1}{2}\begin{pmatrix}
\bm{x}\\
\bm{y}
\end{pmatrix}^{\top}\bm{M}_{tt_0}\begin{pmatrix}
\bm{x}\\
\bm{y}
\end{pmatrix}\right),
\label{kappaintermedform}	
\end{align}
is to compute the pre-factor $c(t,t_0)$. Unique identification of $c(t,t_0)$ serves the dual purpose of uniquely determining the kernel as well as verifying the postulated form \eqref{kappaExpDist}.

As was the case in the last example in Sec. \ref{sec:distfromocp}, here $\kappa$ is not a transition probability kernel and $\int \kappa \neq 1$. Thus, the pre-factor $c(t,t_0)$ cannot be determined from normalization. To find $c(t,t_0)$, the idea now is to substitute \eqref{kappaExpDist} with \eqref{HalfDistSquaredABQ} in \eqref{ReactionAdvectionDiffusionPDEIVP} and invoke the initial condition. 

Specifically, here $$L\kappa\equiv-\langle\nabla_{\bm{x}},\kappa\bm{A}_{t}\bm{x}\rangle + \bm{B}_{t}\bm{B}_{t}^{\top}\Delta_{\bm{x}}\kappa,$$ so $L$ is an advection-diffusion operator. Also, recall that the reaction rate $q(\bm{x})\equiv\frac{1}{2}\bm{x}^{\top}\bm{Q}_{\tau}\bm{x}$ with $\bm{Q}_{\tau}$ satisfying Assumption \textbf{A2}. So \eqref{ReactionAdvectionDiffusionPDEIVP} becomes an \emph{reaction-advection-diffusion PDE initial value problem}
\begin{subequations}
\begin{align}
&\partial_{t}\kappa = -\langle\nabla_{\bm{x}},\kappa\bm{A}_{t}\bm{x}\rangle + \langle\bm{B}_{t}\bm{B}_{t}^{\top},\nabla_{\bm{x}}^{2}\kappa\rangle - \frac{1}{2}\bm{x}^{\top}\bm{Q}_t\bm{x}\kappa, \label{AdvectionReactionDiffusionPDE}\\
&\kappa(t_{0},\bm{x},t_{0},\bm{y})=\delta(\bm{x}-\bm{y}), \label{AdvectionReactionDiffusionIC}
\end{align}
\label{LTVAdvectionReactionDiffusionPDEIVP}
\end{subequations}
where $\nabla_{\bm{x}}^{2}$ denotes the Hessian w.r.t. $\bm{x}$.

We now state our main result in the following theorem (proof in Appendix \ref{AppProofOfThm:kernelLQ}).
\begin{theorem}[Kernel for LQ non-Gaussian steering]\label{thm:kernelLQ}
For $\tau\in\left[t_0,t\right]$, consider assumptions \textbf{A1}-\textbf{A2}. Let
\begin{subequations}
\begin{align}
\!\!\theta(\tau)&:=\tr\:\bm{A}_{\tau} + \langle \bm{B}_{\tau}\bm{B}_{\tau}^{\top},\hat{\bm{\Phi}}_{{\tau t_0}}^{\top}\hat{\bm{\Gamma}}_{{\tau t_0}}^{-1}\hat{\bm{\Phi}}_{{\tau t_0}} + \bm{\Pi}(t_0, \bm{0}, \tau)\rangle,
\label{deftheta}\\
&=\tr\left(\bm{A}_{\tau} + \bm{B}_{\tau}\bm{B}_{\tau}^{\top}\bm{M}_{11}(\tau,t_{0})\right),
\label{defthetaM11}
\end{align}    
\end{subequations}
\begin{align}
a := \left(2\pi\right)^{-n/2}\displaystyle\lim_{t\downarrow t_0}\det\left(\bm{M}_{11}^{1/2}(t,t_0)\exp\int_{t_0}^{t}\left(\bm{A}_{\tau} +\right.\right.\nonumber\\
\left.\left.\quad\bm{B}_{\tau}\bm{B}_{\tau}^{\top}\bm{M}_{11}(\tau,t_{0})\right)\differential\tau\vphantom{\det\left(\bm{M}_{11}^{1/2}(t,t_0)\exp\int_{t_0}^{t}\left(\bm{A}_{\tau} +\right.\right.}\right),
\label{defa}
\end{align}
where $\hat{\bm{\Phi}}_{\tau t_0},\hat{\bm{\Gamma}}_{\tau t_0},\bm{\Pi}(t_0, \bm{0}, \tau)$ are as in Sec. \ref{subsec:gettingdist}, and $\exp$ denotes the matrix exponential. The matrix $\bm{M}_{11}(\tau,t_{0})$ appearing in \eqref{defthetaM11} and \eqref{defa} is the $(1,1)$ block of $\bm{M}_{\tau t_0}$ given by \eqref{defMtt0}.
Then the Markov kernel $\kappa$ solving the reaction-advection-diffusion PDE initial value problem \eqref{LTVAdvectionReactionDiffusionPDEIVP} is
\begin{align}
\kappa\left(t_0,\bm{x},t,\bm{y}\right)
=& a\exp\left(-\int_{t_0}^{t}\theta(s)\differential s\right)\nonumber\\
&\times\exp\!\left(\!-\frac{1}{2}\!\begin{pmatrix}
\bm{x}\\ \bm{y}
\end{pmatrix}^{\!\!\top}\!\!\bm{M}_{tt_{0}}\!\!\begin{pmatrix}
\bm{x}\\ \bm{y}
\end{pmatrix}\!\right),
\label{KernelInTermsOfM}
\end{align}
where $\bm{M}_{tt_{0}}\succ\bm{0}$ is given by \eqref{defMtt0}.
\end{theorem}
The kernel \eqref{KernelInTermsOfM} significantly extends the results in \cite{teter2024schr,teter2024weyl} to the generic LQ case. The corollary next (proof in Appendix \ref{AppProofOfCorollarySplCase}) recovers the previously known special case.

\begin{corollary}\label{Corollary:SplCaseOfLQkernel}[Kernel in \cite{teter2024schr,teter2024weyl} as special case of \eqref{KernelInTermsOfM}]
When $\bm{Q}_{\tau}=2\bm{D}\succ\bm{0}$ (constant positive diagonal matrix) and $(\bm{A}_{\tau},\hat{\bm{B}}_{\tau})\equiv(\bm{0},\sqrt{2}\bm{I})$ $\forall \tau\in[t_0,t]$, the kernel \eqref{KernelInTermsOfM} reduces to that in \cite[Eq. (43)]{teter2024weyl} or that in \cite[Eq. (A.22)]{teter2024schr}.
\end{corollary}

\begin{remark}\label{Remark:LimitIna0}
Using assumptions \textbf{A1}-\textbf{A2} and the standard $(\varepsilon,\delta)$ definition of limit, it is not too difficult to show that the upper limit in \eqref{defa} exists. In special cases such as Corollary \ref{Corollary:SplCaseOfLQkernel} above, the limit in \eqref{defa} can be evaluated in closed form, as shown in Appendix \ref{AppProofOfCorollarySplCase}. Importantly, the limit in \eqref{defa} cannot be distributed to the determinant and the exponential factors; see the computation in \eqref{aSplCase}. 
\end{remark}

The following properties are immediate from the expression of the kernel \eqref{KernelInTermsOfM}:
\begin{itemize}
\item (spatial symmetry) $\kappa\left(t_0,\bm{x},t,\bm{y}\right) = \kappa\left(t_0,\bm{y},t,\bm{x}\right)$ $\forall\bm{x},\bm{y}\in\mathbb{R}^{n}$,

\item (positivity) $\kappa\left(t_0,\bm{x},t,\bm{y}\right)>0$ since from \eqref{defa}, $a>0$.
\end{itemize}

We next discuss how the derived kernel helps in solving the LQ SBP \eqref{LinearSBPstatecost}. 


\subsection{Using $\kappa$ to solve LQ SBP with non-Gaussian $\mu_0,\mu_1$}\label{subsec:FromkappatoBridge}
Using the derived kernel \eqref{KernelInTermsOfM}, we define the integral transforms
\begin{subequations}
\begin{align}
&\widehat{\varphi}_{0} \mapsto \mathcal{T}_{01}\widehat{\varphi}_{0}:=\int_{\mathbb{R}^{n}}\kappa\left(t_0,\bm{x},t_{1},\bm{y}\right)\widehat{\varphi}_{0}(\bm{y}) \differential\bm{y},\label{IntegralMapForward}\\
&\varphi_{1} \mapsto \mathcal{T}_{10}\varphi_{1}:=\int_{\mathbb{R}^{n}}\kappa\left(t_1,\bm{x},t_{0},\bm{y}\right)\varphi_{1}(\bm{y}) \differential\bm{y},\label{IntegralMapBackward}
\end{align}
\label{DefIntegralOperator}   
\end{subequations}
for measurable $\widehat{\varphi}_{0}(\cdot),\varphi_{1}(\cdot)$. These integral transforms help solve the LQ SBP \eqref{LinearSBPstatecost} as follows.

Suppose that the endpoint measures $\mu_0,\mu_1$ in \eqref{LinearSBPstatecostEndpointConstr} are absolutely continuous with respective probability density functions (PDFs) $\rho_0,\rho_1$. Standard computation shows that the necessary conditions for optimality for \eqref{LinearSBPstatecost} yield a pair of second order coupled nonlinear PDEs: 
\begin{subequations}
\begin{align}
&\partial_{t} \rho_{\mathrm{opt}}^{\bm{u}}+\nabla_{\boldsymbol{x}} \cdot\left(\rho_{\mathrm{opt}}^{\bm{u}}\left(\boldsymbol{A}_{t}\bm{x}+\boldsymbol{B}_{t} \boldsymbol{B}_{t}^{\top} \nabla_{\boldsymbol{x}} S\right)\right)\nonumber\\
&\hspace*{1.8in}=\langle\bm{B}_{t}\bm{B}_{t}^{\top},\nabla_{\bm{x}}^{2} \rho_{\mathrm{opt}}^{\bm{u}}\rangle,\label{ControlledFPK}\\
&\partial_{t} S +\left\langle\nabla_{\boldsymbol{x}} S, \boldsymbol{A}_{t}\bm{x}\right\rangle+\frac{1}{2}\left\langle\nabla_{\boldsymbol{x}} S, \boldsymbol{B}_{t} \boldsymbol{B}_{t}^{\top} \nabla_{\boldsymbol{x}} S\right\rangle\nonumber\\
&\hspace*{1in}+\left\langle\bm{B}_{t}\bm{B}_{t}^{\top}, \nabla_{\bm{x}}^{2} S\right\rangle=\frac{1}{2}\bm{x}^{\top}\bm{Q}_{t}\bm{x},\label{HJBlike}\\
&\rho_{\mathrm{opt}}^{\bm{u}}\left(t=t_0, \cdot\right)=\rho_0(\cdot), \quad \rho_{\mathrm{opt}}^{\bm{u}}\left(t=t_1, \cdot\right)=\rho_1(\cdot),\label{EndpointBC}
\end{align}
\label{FOOC}
\end{subequations}
in unknown primal-dual pair $\left(\rho_{\mathrm{opt}}^{\bm{u}}(t,\bm{x}),S(t,\bm{x})\right)$, i.e., the optimally controlled joint state PDF and the value function.

The Hopf-Cole transform \cite{hopf1950partial,cole1951quasi} given by
\begin{align}
\left(\rho_{\mathrm{opt}}^{\bm{u}},S\right) \mapsto \left(\widehat{\varphi},\varphi\right):=\left(\rho_{\mathrm{opt}}^{\bm{u}}\exp\left(-S\right),\exp S\right),
\label{HopfColeTransform}    
\end{align}
recasts \eqref{FOOC} to a pair of boundary-coupled linear PDEs:
\begin{subequations}
\begin{align}
&\!\!\partial_{t}\widehat{\varphi} = -\langle\nabla_{\bm{x}},\widehat{\varphi}\bm{A}_{t}\bm{x}\rangle + \langle\bm{B}_{t}\bm{B}_{t}^{\top},\nabla_{\bm{x}}^{2}\widehat{\varphi}\rangle - \frac{1}{2}\bm{x}^{\top}\bm{Q}_t\bm{x}\widehat{\varphi},\label{phihatPDE}\\
&\!\!\partial_{t}\varphi = -\langle\nabla_{\bm{x}}\varphi,\bm{A}_{t}\bm{x}\rangle - \langle\bm{B}_{t}\bm{B}_{t}^{\top},\nabla_{\bm{x}}^{2}\varphi\rangle + \frac{1}{2}\bm{x}^{\top}\bm{Q}_t\bm{x}\varphi, \label{phiPDE}\\
&\!\!\widehat{\varphi}(t_0,\cdot)\varphi(t_0,\cdot)=\rho_0(\cdot), \quad \widehat{\varphi}(t_1,\cdot)\varphi(t_1,\cdot)=\rho_1(\cdot).\label{phihatphiBC}   
\end{align}
\label{phihatphiPDE}    
\end{subequations}
We note that the PDE \eqref{phihatPDE} is precisely the forward reaction-advection-diffusion PDE \eqref{AdvectionReactionDiffusionPDE}, and the PDE \eqref{phiPDE} is the associated backward reaction-advection-diffusion PDE. 

The solution of \eqref{phihatphiPDE} recovers the solution of \eqref{LinearSBPstatecost}. Specifically, $\forall t\in[t_0,t_1]$, the optimally controlled joint state PDF $\rho_{\mathrm{opt}}^{\bm{u}}(t,\bm{x}) = \widehat{\varphi}(t,\bm{x})\varphi(t,\bm{x})$, and the optimal control $\bm{u}_{\mathrm{opt}}(t,\bm{x}) = \bm{B}_{t}^{\top}\nabla_{\bm{x}}S(t,\bm{x}) =  \bm{B}_{t}^{\top}\nabla_{\bm{x}}\log\varphi(t,\bm{x})$.

For $k\in\{0,1\}$, letting $\widehat{\varphi}_{k}(\cdot):= \widehat{\varphi}(t_k,\cdot)$, $\varphi_{k}(\cdot):=\varphi(t_k,\cdot)$, the system \eqref{phihatphiPDE} can be solved by the dynamic Sinkhorn recursion:
\begin{align}
\widehat{\varphi}_{0} ~\xrightarrow{\mathcal{T}_{01}}~ \widehat{\varphi}_{1} ~\xrightarrow{\rho_1/\widehat{\varphi}_{1}}~ \varphi_{1} ~\xrightarrow{\mathcal{T}_{10}}~ \varphi_{0} ~\xrightarrow{\rho_0/\varphi_{0}}~ \left(\widehat{\varphi}_{0}\right)_{\text{next}},
\label{DynamicSinkhornRecursion}    
\end{align}
involving the integral transforms \eqref{DefIntegralOperator}. The recursion \eqref{DynamicSinkhornRecursion} is known \cite{chen2016entropic,marino2020optimal} to be contractive w.r.t. Hilbert's projective metric \cite{bushell1973hilbert,lemmens2014birkhoff}. For different variants of the dynamic Sinkhorn recursions in the SBP context, we refer the readers to \cite{caluya2021wasserstein,caluya2021reflected,teter2023contraction,teter2025probabilistic}. For discussions on \eqref{DynamicSinkhornRecursion} in the special case $\left(\bm{A}_{t},\bm{B}_{t}\right)\equiv\left(\bm{0},\bm{I}\right)$, see \cite[Sec. 3.2]{teter2024schr}.

From computational point of view, having a handle for the kernel \eqref{KernelInTermsOfM} helps us to apply $\mathcal{T}_{01},\mathcal{T}_{10}$ in \eqref{DefIntegralOperator} during each pass of the recursion \eqref{DynamicSinkhornRecursion}. This in turn helps to solve the LQ SBP \eqref{LinearSBPstatecost} with arbitrary \emph{non-Gaussian} $\mu_0,\mu_1$ having finite second moments. The solution for \eqref{LinearSBPstatecost} when both $\mu_0,\mu_1$ are \emph{Gaussians}, appeared in \cite{chenPartIII} and did not need to derive the kernel. However, the techniques and results in \cite{chenPartIII} are difficult to generalize for non-Gaussian $\mu_0,\mu_1$. This is what motivated our derivation for the kernel.

\begin{remark}
If the kernel $\kappa$ in \eqref{DefIntegralOperator} were not available, an alternative way to implement $\mathcal{T}_{01},\mathcal{T}_{10}$ in \eqref{DynamicSinkhornRecursion} is to apply the Feynman-Kac path integral \cite[Ch. 8.2]{oksendal2013stochastic}, \cite[Ch. 3.3]{yong2012stochastic} resulting in randomized numerical approximations for $\widehat{\varphi}_{1},\varphi_{0}$. For instance, the Feynman-Kac path integral representation for \eqref{IntegralMapBackward} is
\begin{align}
&\mathcal{T}_{10}\varphi_{1}(\bm{x})\nonumber\\
=&\mathbb{E}\left[\varphi_{1}\left(\bm{x}_{t_{1}}\right)\exp\left(-\int_{t}^{t_1}\frac{1}{2}\bm{x}_{\tau}^{\top}\bm{Q}_{\tau}\bm{x}_{\tau}\differential\tau\right)\mid \bm{x}_{t}=\bm{x}\right],
\end{align}
and the conditional expectation can be approximated by Monte Carlo simulation. Having an explicit handle on $\kappa$ removes the need for such randomized approximation of the deterministic functions $\widehat{\varphi}_{1},\varphi_{0}$.
\end{remark}

To further understand the action of the derived kernel, notice that combining \eqref{KernelInTermsOfM} and \eqref{IntegralMapForward} yields
\begin{align}
&\mathcal{T}_{01}\widehat{\varphi}_{0} = a\exp\left(-\int_{t_0}^{t}\theta(s)\differential s\right)\exp\left(-\frac{1}{2}\bm{x}^{\top}\bm{M}_{11}\left(t,t_0\right)\bm{x}\right)\nonumber\\
&\int_{\mathbb{R}^{n}}\!\!\exp\!\left(\!-\frac{1}{2}\bm{y}^{\top}\bm{M}_{22}\left(t,t_{0}\right)\bm{y} + \langle-\bm{M}_{12}^{\top}\left(t,t_0\right)\bm{x},\bm{y} \rangle\!\right)\widehat{\varphi}_{0}(\bm{y})\differential\bm{y},
\label{T01explicit}    
\end{align}
and likewise for $\mathcal{T}_{10}\varphi_{1}$, where $\bm{M}_{11},\bm{M}_{12},\bm{M}_{22}$ refer to the respective blocks of \eqref{defMtt0}. For suitably smooth $\widehat{\varphi}_{0}$, the integral in \eqref{T01explicit} can be evaluated using Lemma \ref{Lemma:GaussianIntegral} in Appendix \ref{DerivationOfcClassicalReactionDiffusion}. We close with an example of such computation.

\begin{example}[Mixture-of-Gaussian $\widehat{\varphi}_{0}$]
For fixed $n_{c}\in\mathbb{N}$, let $\widehat{\varphi}_{0}$ be an $n_{c}$ component conic mixture-of-Gaussians, i.e., $$\widehat{\varphi}_{0}(\bm{y})=\sum_{i=1}^{n_{c}} w_{i}\exp\left(-\frac{1}{2}\left(\bm{y}-\bm{m}_{i}\right)^{\top}\bm{\Sigma}_{i}^{-1}\left(\bm{y}-\bm{m}_{i}\right)\right),$$
where $w_{i}>0$, $\bm{m}_{i}\in\mathbb{R}^{n}$, $\bm{\Sigma}_{i}\succ\bm{0}$ $\forall i\in[n_{c}]$. Then \eqref{T01explicit} gives
\begin{align}
&\mathcal{T}_{01}\widehat{\varphi}_{0} = a\exp\left(-\int_{t_0}^{t}\theta(s)\differential s\right)\nonumber\\
&\times\exp\left(-\frac{1}{2}\bm{x}^{\top}\bm{M}_{11}\left(t,t_0\right)\bm{x}\right)\sum_{i=1}^{n_{c}}w_{i}\exp\left(-\frac{1}{2}\bm{m}_{i}^{\top}\bm{\Sigma}_{i}^{-1}\bm{m}_{i}\right)\nonumber\\
&\int_{\mathbb{R}^{n}}\!\!\exp\!\left(\!-\frac{1}{2}\bm{y}^{\top}\left(\bm{M}_{22}\left(t,t_{0}\right)+\bm{\Sigma}_{i}^{-1}\right)\bm{y} \right.\nonumber\\
&\left.\qquad\qquad+ \langle-\bm{M}_{12}^{\top}\left(t,t_0\right)\bm{x}+\bm{\Sigma}_{i}^{-1}\bm{m}_{i},\bm{y} \rangle\vphantom{-\frac{1}{2}\bm{y}^{\top}\left(\bm{M}_{22}\left(t,t_{0}\right)+\bm{\Sigma}_{i}^{-1}\right)\bm{y}}\!\right)\differential\bm{y}\nonumber\\
&= a\left(2\pi\right)^{n/2}\exp\left(-\int_{t_0}^{t}\theta(s)\differential s\right)\nonumber\\
&\times\exp\left(-\frac{1}{2}\bm{x}^{\top}\bm{M}_{11}\left(t,t_0\right)\bm{x}\right)\sum_{i=1}^{n_{c}}\dfrac{w_{i}}{\sqrt{\det\left(\bm{M}_{22}\left(t,t_{0}\right)+\bm{\Sigma}_{i}^{-1}\right)}}\nonumber\\
&\times\!\exp\bigg\{\!-\frac{1}{2}\bm{m}_{i}^{\top}\bm{\Sigma}_{i}^{-1}\bm{m}_{i}+\frac{1}{2}\left(-\bm{M}_{12}^{\top}\left(t,t_0\right)\bm{x}+\bm{\Sigma}_{i}^{-1}\bm{m}_{i}\right)^{\top}\nonumber\\
&\qquad\quad\left(\bm{M}_{22}\left(t,t_{0}\right)+\bm{\Sigma}_{i}^{-1}\right)\left(-\bm{M}_{12}^{\top}\left(t,t_0\right)\bm{x}+\bm{\Sigma}_{i}^{-1}\bm{m}_{i}\right)\!\bigg\},
\end{align}
where the last equality uses Lemma \ref{Lemma:GaussianIntegral}.
\end{example}

\section{Concluding Remarks}\label{sec:conclusions}
We derived the Markov kernel for an LTV It\^{o} diffusion in the presence of killing of probability mass at a rate that is convex quadratic with time-varying weight matrix. The resulting kernel is parameterized by the LTV matrix pair $(\bm{A}_{t},\bm{B}_{t})$ and the killing weight $\bm{Q}_{t}\succeq\bm{0}$, in terms of the solution of a Riccati matrix ODE initial value problem. 

The derived kernel has relevance in stochastic control: it is the Green's function for a reaction-advection-diffusion PDE that appears in solving the generic linear quadratic Schr\"{o}dinger bridge problem. This problem concerns with steering a given state distribution to another over a given finite horizon subject to controlled LTV diffusion while minimizing a cost-to-go that is quadratic in state and control. The solution for this problem has appeared in prior literature for the case when the endpoint distributions are Gaussians. It is also understood that for non-Gaussian endpoints, the problem can be solved via dynamic Sinkhorn recursions, which however, require solving initial value problems involving the aforesaid reaction-advection-diffusion PDE within each epoch of the recursion, with updated initial conditions. By deriving the corresponding Green's function, our results facilitate this computation.

The results here also generalize our prior works where the Markov kernel was derived for the special case $(\bm{A}_{t},\bm{B}_{t})=\left(\bm{0},\bm{I}\right)$ using generalized Hermite polynomials and Weyl calculus. However, for generic $(\bm{A}_{t},\bm{B}_{t})$ and $\bm{Q}_{t}\succeq\bm{0}$, those techniques become unwieldy. To overcome this technical challenge, we pursued a new line of attack by postulating the structure of the Markov kernel in terms of a distance function induced by an underlying deterministic optimal control problem. Using this new technique, we demonstrated that both new and existing results can be recovered in a conceptually transparent manner even when the underlying Markov kernel is not a transition probability, as is the case here. These interconnections between the Markov kernels, distances and optimal control, should be of independent interest.   

\appendix

\subsection{Derivation of \eqref{PrefactorKernelSICON}}\label{DerivationOfcClassicalReactionDiffusion}
Consider $\kappa$ as in \eqref{kappaExpDist} with $\frac{1}{2}\dist^{2}_{tt_{0}}$ given by \eqref{halfdistsquaredSICONkernel}. For finding the pre-factor $c(t,t_{0})$, we substitute \eqref{kappaExpDist} with \eqref{halfdistsquaredSICONkernel} in the reaction-diffusion PDE in \eqref{ReactionAdvectionDiffusionPDEIVP}, i.e., in \begin{align}
\partial_{t}\kappa &= \left(\Delta_{\bm{x}} -  \frac{1}{2}\bm{x}^{\top}2\bm{D}\bm{x}\right)\kappa\nonumber\\
&= \left(\displaystyle\sum_{i=1}^{n}\partial^{2}_{x_{ii}} -  \dfrac{\omega_i^2}{4}x_{ii}^{2}\right)\kappa,
\label{ReactionDiffusionClassicalSpectral}
\end{align}
since $\omega_{i}^{2}=4D_{ii}$ for all $i=1,\hdots,n$. This substitution in \eqref{ReactionDiffusionClassicalSpectral}, after algebraic simplification, results in the ODE 
\begin{align}
\dot{c} = c\left(-\dfrac{1}{2}\displaystyle\sum_{i=1}^{n}\omega_{i}\coth\left(\omega_{i}(t-t_{0})\right)\right),
\label{cdotODE}
\end{align}
where $\dot{c}$ denotes the derivative of $c(t,t_0)$ w.r.t. $t$. 

Integrating \eqref{cdotODE} yields
$$\log c = \log a -\dfrac{1}{2}\displaystyle\sum_{i=1}^{n}\log\sinh(\omega_i (t-t_0)),$$
where $\log a$ is the numerical constant of integration. Thus,
\begin{align}
c(t,t_0) = a\displaystyle\prod_{i=1}^{n}\dfrac{1}{\left(\sinh(\omega_i (t-t_0))\right)^{1/2}}.
\label{cwitha}    
\end{align}
Combining \eqref{cwitha} with \eqref{kappaExpDist} and \eqref{halfdistsquaredSICONkernel}, we obtain
\begin{align}
&\kappa(t_{0},\bm{x},t,\bm{y}) = a\displaystyle\prod_{i=1}^{n}\dfrac{1}{\sqrt{\sinh(\omega_i (t-t_0))}}\nonumber\\
&\times\exp\left(-\dfrac{1}{2}\dfrac{\omega_{i}\left(x_i^2 + y_i^2\right)\cosh\left(\omega_{i}(t-t_{0})\right) - 2\omega_i x_i y_i}{2\sinh\left(\omega_{i}(t-t_{0})\right)}\right),
\label{kappawitha}    
\end{align}
for all $0\leq t_{0}\leq t < \infty$. All that remains is to find the constant $a$.

Letting $\tau:= t-t_{0}$, and invoking the initial condition in \eqref{ReactionAdvectionDiffusionPDEIVP} for \eqref{kappawitha}, we have
\begin{align}
&\delta(\bm{x}-\bm{y}) = a\displaystyle\lim_{\tau\downarrow 0} \prod_{i=1}^{n}\dfrac{1}{\sqrt{\sinh(\omega_i \tau)}}\nonumber\\
&\times \exp\left(-\dfrac{1}{2}\dfrac{\omega_{i}\left(x_i^2 + y_i^2\right)\cosh\left(\omega_{i}\tau\right) - 2\omega_i x_i y_i}{2\sinh\left(\omega_{i}\tau\right)}\right)
\label{FindingcClassicalReactionDiffusionIC}.    
\end{align}
Integrating both sides of \eqref{FindingcClassicalReactionDiffusionIC} w.r.t. $\bm{x}$ over $\mathbb{R}^{n}$ gives 
\begin{align}
&1 = a\left(\displaystyle\lim_{\tau\downarrow 0}\prod_{i=1}^{n}\dfrac{1}{\sqrt{\sinh(\omega_i \tau)}}\right)\times\nonumber\\
&\!\int_{\mathbb{R}^{n}}\!\displaystyle\lim_{\tau\downarrow 0}\prod_{i=1}^{n}\!\exp\!\left(\!-\dfrac{1}{2}\dfrac{\omega_{i}\left(x_i^2 + y_i^2\right)\cosh\left(\omega_{i}\tau\right) - 2\omega_i x_i y_i}{2\sinh\left(\omega_{i}\tau\right)}\!\right)\!\differential\bm{x},
\label{Findinga:start}
\end{align}
wherein the LHS used the shift property of the Dirac delta: $\int_{\mathbb{R}^{n}}\delta(\bm{x}-\bm{y})\cdot 1\cdot\differential\bm{x} = 1$.

Since  
\begin{align*}
&\dfrac{\omega_{i}\left(x_i^2 + y_i^2\right)\cosh\left(\omega_{i}\tau\right) - 2\omega_i x_i y_i}{2\sinh\left(\omega_{i}\tau\right)}\\
&= \begin{pmatrix}
\sqrt{\omega}_ix_i \\ \sqrt{\omega}_i y_i
\end{pmatrix}^{\top} \begin{pmatrix}
\coth(\omega_i\tau) & -{\mathrm{csch}}(\omega_i\tau)\\
-{\mathrm{csch}}(\omega_i\tau) & \coth(\omega_i\tau)
\end{pmatrix}
\begin{pmatrix}
\sqrt{\omega}_ix_i \\ \sqrt{\omega}_i y_i
\end{pmatrix}
\end{align*}
is a positive definite quadratic form, we have
$$0\leq \exp\left(-\dfrac{1}{2}\dfrac{\omega_{i}\left(x_i^2 + y_i^2\right)\cosh\left(\omega_{i}\tau\right) - 2\omega_i x_i y_i}{2\sinh\left(\omega_{i}\tau\right)}\right)\leq 1.$$
Hence, using the dominated convergence theorem \cite[Thm. 1.13]{stein2009real}, we exchange the limit and integral in \eqref{Findinga:start} to yield
\begin{align}
1 &= a\displaystyle\lim_{\tau\downarrow 0}\prod_{i=1}^{n}\dfrac{\exp\left(-\dfrac{1}{2}\dfrac{\omega_i y_i^2 \cosh(\omega_i\tau)}{2\sinh(\omega_i\tau)}\right)}{\sqrt{\sinh(\omega_i\tau)}}\nonumber\\
&\times\left(\int_{-\infty}^{\infty}\exp\left(-\dfrac{1}{2}\dfrac{\omega_i x_i^2 \cosh(\omega_i\tau) - 2\omega_i x_i y_i}{2\sinh(\omega_i\tau)}\right)\differential x_{i}\right).
\label{Findinga:almostend}    
\end{align}
We need the following auxiliary lemma.
\begin{lemma}[Central identity of quantum field theory] \cite[p. 2]{zinn2021quantum}, \cite[p. 15]{zee2010quantum}\label{Lemma:GaussianIntegral} For fixed $n\times n$ matrix $\bm{A}\succ\bm{0}$ and suitably smooth $f:\mathbb{R}^{n}\mapsto\mathbb{R}$, 
\begin{align}
&\int_{\mathbb{R}^n} \exp \left(-\frac{1}{2} \boldsymbol{x}^{\top} \boldsymbol{A} \boldsymbol{x}\right) f(\boldsymbol{x}) \mathrm{d} \boldsymbol{x}\nonumber\\
&=\left.\sqrt{\frac{(2 \pi)^n}{\operatorname{det}(\boldsymbol{A})}} \exp \left(\frac{1}{2} \nabla_{\boldsymbol{x}}^{\top} \boldsymbol{A}^{-1} \nabla_{\boldsymbol{x}}\right) f(\boldsymbol{x})\right|_{\boldsymbol{x}=\mathbf{0}},
\label{CentralIdentityQFT}   
\end{align}
where the exponential of differential operator is understood as a power series.
As a special case, for $\bm{b}\in\mathbb{R}^{n}$, we have
$\int_{\mathbb{R}^{n}}\exp\left(-\frac{1}{2}\bm{x}^{\top}\bm{A}\bm{x}+\langle\bm{b},\bm{x}\rangle\right)\differential\bm{x}
= \sqrt{\frac{(2\pi)^{n}}{\det\bm{A}}}\exp\left(\frac{1}{2}\bm{b}^{\top}\bm{A}^{-1}\bm{b}\right)$. In particular, for fixed scalars $a>0, b\in\mathbb{R}$, we have $$\int_{-\infty}^{\infty}\exp\left(-\frac{1}{2}a x^2 + b x\right)\differential x = \sqrt{\frac{2\pi}{a}}\exp\left(\frac{1}{2}\frac{b^2}{a}\right).$$
\end{lemma}

Using Lemma \ref{Lemma:GaussianIntegral}, the integral in \eqref{Findinga:almostend} evaluates to
$$\sqrt{\dfrac{4\pi\sinh(\omega_i\tau)}{\omega_i\cosh(\omega_i\tau)}}\exp\left(\dfrac{\omega_i y_i^2}{4\sinh(\omega_i\tau)\cosh(\omega_i\tau)}\right).$$
Substituting back this value in \eqref{Findinga:almostend}, followed by simplification via identity: $1 - \cosh^{2}(\cdot) = -\sinh^{2}(\cdot)$, gives
\begin{align}
1 &= a\displaystyle\lim_{\tau\downarrow 0}\prod_{i=1}^{n}\sqrt{\dfrac{4\pi}{\omega_i\cosh(\omega_i\tau)}}\exp\left(-\dfrac{\omega_i y_i^2\sinh(\omega_i\tau)}{4}\right)\nonumber\\
&= a\prod_{i=1}^{n}\sqrt{\dfrac{4\pi}{\omega_i}}.
\label{Findinga:end}
\end{align}
Therefore, $a = \prod_{i=1}^{n}\sqrt{\frac{\omega_i}{4\pi}}$, which upon substitution in \eqref{cwitha}, yields \eqref{PrefactorKernelSICON}.


\subsection{Proof of Proposition \ref{prop:SpecialCaseAzeroBidentity}}\label{AppProofOfProp:SpecialCaseAzeroBidentity}
As is well-known \cite[p. 156]{brockett1970finite}, \cite{shayman1986phase}, for any $\tau\in[t_0,t]$, the solution of the Riccati ODE initial value problem \eqref{RiccatiIVP} admits linear fractional representation
\begin{align}
&\bm{\Pi}(\tau,\bm{K}_{1},t)= \nonumber\\
&\left(\bm{\Psi}_{21}\left(t,\tau\right) + \bm{\Psi}_{22}\left(t,\tau\right)\bm{K}_{1}\right)\left(\bm{\Psi}_{11}\left(t,\tau\right) + \bm{\Psi}_{12}\left(t,\tau\right)\bm{K}_{1}\right)^{-1},
\label{PiLFT}    
\end{align}
where 
\begin{align}
\bm{\Psi}_{t\tau} := \begin{bmatrix}
\bm{\Psi}_{11}\left(t,\tau\right) & \bm{\Psi}_{12}\left(t,\tau\right)\\
\bm{\Psi}_{21}\left(t,\tau\right) & \bm{\Psi}_{22}\left(t,\tau\right)
\end{bmatrix}\in\mathbb{R}^{2n\times 2n}
\label{defPsiBlock}    
\end{align}
is the state transition matrix of the linear Hamiltonian matrix ODE
\begin{align}
\begin{pmatrix}
\dot{\bm{X}}_{\tau}\\
\dot{\bm{\Lambda}}_{\tau}
\end{pmatrix} = \begin{bmatrix}
\bm{A}_{\tau} & -\hat{\bm{B}}_{\tau}\hat{\bm{B}}_{\tau}^{\top}\\
-\bm{Q}_{\tau} & -\bm{A}_{\tau}^{\top}
\end{bmatrix}\begin{pmatrix}
\bm{X}_{\tau}\\
\bm{\Lambda}_{\tau}
\end{pmatrix}, \quad \bm{X}_{\tau},\bm{\Lambda}_{\tau}\in\mathbb{R}^{n\times n}.
\label{HamiltonianMatrixODE}
\end{align}
For the special case in hand, the coefficient matrix in \eqref{HamiltonianMatrixODE} equals 
$$\begin{bmatrix}
\bm{0} & -2\bm{I}\\
-2\bm{D} & \bm{0}
\end{bmatrix},$$
and its (backward in time) state transition matrix becomes
\begin{align}
&\bm{\Psi}_{t\tau} = \exp\left(\begin{bmatrix}
 \bm{0} & 2\bm{I}\\
2\bm{D} & \bm{0}   
\end{bmatrix}(t-\tau)\right) \nonumber\\
&=\!\! {\small{\begin{bmatrix}
\cosh\left(2\sqrt{\bm{D}}(t-\tau)\right) & \left(\sqrt{\bm{D}}\right)^{\!-1}\!\!\sinh\left(2\sqrt{\bm{D}}(t-\tau)\right)\\
\sqrt{\bm{D}}\sinh\left(2\sqrt{\bm{D}}(t-\tau)\right) & \cosh\left(2\sqrt{\bm{D}}(t-\tau)\right)
\end{bmatrix}}},
\label{Psittau}
\end{align}
where all hyperbolic functions act element-wise. To see the last equality, it suffices to note that for $2\times 2$ matrix $\begin{bmatrix}
0 & b\\
c & 0
\end{bmatrix}$ with $b,c>0$, we have
\begin{align*}
&\exp\left(\begin{bmatrix}
0 & b\\
c & 0
\end{bmatrix}(t-\tau)\right)= \bm{I}\sum_{k=0}^{\infty}\frac{\left(\sqrt{bc}(t-\tau)\right)^{2k}}{(2k)!}\\
&\qquad\qquad\qquad + \begin{bmatrix}
0 & \sqrt{b/c}\\
\sqrt{c/b} & 0
\end{bmatrix}\sum_{k=0}^{\infty}\frac{\left(\sqrt{bc}(t-\tau)\right)^{2k+1}}{(2k+1)!}\\
&=\begin{bmatrix}
\cosh(\sqrt{bc}(t-\tau)) & \sqrt{b/c}\:\sinh(\sqrt{bc}(t-\tau))\\
\sqrt{c/b}\:\sinh(\sqrt{bc}(t-\tau)) & \cosh(\sqrt{bc}(t-\tau))
\end{bmatrix}.
\end{align*}

From \eqref{PiLFT} and \eqref{Psittau}, we then get
\begin{align}
\bm{\Pi}(\tau,\bm{0},t) &= \bm{\Psi}_{21}\left(t,\tau\right)\left(\bm{\Psi}_{11}\left(t,\tau\right)\right)^{-1} \nonumber\\
&= \sqrt{\bm{D}}\tanh\left(2\sqrt{\bm{D}}(t-\tau)\right).
\label{Pisplcase}
\end{align}

Following \eqref{defAhat}, in this case, we also have 
\begin{align}
\hat{\bm{A}}_{\tau} = -2\bm{\Pi}(\tau,\bm{0},t) = -2\sqrt{\bm{D}}\tanh\left(2\sqrt{\bm{D}}(t-\tau)\right)
\label{Ahattausplcase}    
\end{align}
for all $\tau\in[t_0,t]$. From \eqref{Ahattausplcase}, $\hat{\bm{A}}_{\tau}$ and $\int_{t_0}^{\tau}\hat{\bm{A}}_{\sigma}\differential\sigma$ commute. Hence the state transition matrix for the coefficient matrix $\hat{\bm{A}}_{\tau}$ equals
\begin{align}
\hat{\bm{\Phi}}_{t\tau} &= \exp\left(-2\sqrt{\bm{D}}\int_{\tau}^{t}\tanh\left(2\sqrt{\bm{D}}(t-\sigma)\right)\differential\sigma\right)\nonumber\\
&= \exp\left(-2\sqrt{\bm{D}}\int_{0}^{t-\tau}\tanh\left(2\sqrt{\bm{D}}s\right)\differential s\right)\nonumber\\
&=\exp\left(-\log\cosh\left(2\sqrt{\bm{D}}(t-\tau)\right)\right)\nonumber\\
&= {\mathrm{sech}}\left(2\sqrt{\bm{D}}(t-\tau)\right). 
\label{SpecialCasePhihat}    
\end{align}
Thus 
\begin{align*}
\hat{\bm{\Gamma}}_{t\tau} =\int_{\tau}^{t}\hat{\bm{\Phi}}_{t\sigma}\hat{\bm{B}}_{\sigma}\hat{\bm{B}}_{\sigma}^{\top}\hat{\bm{\Phi}}_{t\sigma}^{\top}\differential\sigma &= 2\int_{\tau}^{t}{\mathrm{sech}}^{2}(2\sqrt{\bm{D}}(t-\sigma))\differential\sigma\\
&= \bm{D}^{-\frac{1}{2}}\tanh(2\sqrt{\bm{D}}(t-\tau)),    
\end{align*} 
so,
\begin{align}
\hat{\bm{\Gamma}}_{t\tau}^{-1} = \sqrt{\bm{D}}\coth(2\sqrt{\bm{D}}(t-\tau)).
\label{SpecialCaseInvGammahat}    
\end{align}
Therefore, \eqref{defMtt0} simplifies to
\begin{align}
\bm{M}_{tt_0}\!\!=\!\!{\small{\begin{bmatrix}
\!\sqrt{\bm{D}}\coth\left(2\sqrt{\bm{D}}(t-t_0)\right) & -\sqrt{\bm{D}}{\mathrm{csch}}\left(2\sqrt{\bm{D}}(t-t_0)\right)\!\!\\
\!-\sqrt{\bm{D}}{\mathrm{csch}}\left(2\sqrt{\bm{D}}(t-t_0)\right) & \sqrt{\bm{D}}\coth\left(2\sqrt{\bm{D}}(t-t_0)\right)\!
\!\end{bmatrix}}}
\label{Nsplcase}    
\end{align}
which is exactly formula (4.10) in \cite{teter2024schr} derived via completely different computation.

Substituting \eqref{Nsplcase} in \eqref{HalfDistSquaredABQ}, and recalling that $\omega_i := 2\sqrt{D_{ii}}$ for all $i=1,\hdots, n$, we recover \eqref{halfdistsquaredSICONkernel}. \qed


\subsection{Proof of Theorem \ref{thm:kernelLQ}}\label{AppProofOfThm:kernelLQ}
The proof is divided into three steps.

\noindent{\textbf{Step 1: Determining $c(t,t_0)$ up to a constant $a$.}}\\
\noindent Combining \eqref{kappaExpDist} with \eqref{HalfDistSquaredABQ}, we get
\begin{align}
\log\kappa = \log c(t,t_0) - \frac{1}{2}\begin{pmatrix}
\bm{x}\\
\bm{y}
\end{pmatrix}^{\top}\bm{M}_{tt_0}\begin{pmatrix}
\bm{x}\\
\bm{y}
\end{pmatrix}.
\label{logofkernel}
\end{align}
Applying $\partial_t$ to both sides of \eqref{logofkernel} and rearranging, gives
\begin{align}
\partial_t\kappa = \kappa\left(\frac{\dot{c}}{c} - \frac{1}{2}\begin{pmatrix}
\bm{x}\\
\bm{y}
\end{pmatrix}^{\top}\dot{\bm{M}}_{tt_0}\begin{pmatrix}
\bm{x}\\
\bm{y}
\end{pmatrix}\right),
\label{LHSadvectionreactiondiffusionPDE}
\end{align}
where dot denotes derivative w.r.t. $t$. 

Applying $\nabla_{\bm{x}}$ to both sides of \eqref{logofkernel}, and using \eqref{defMtt0}, we find
\begin{align}
\nabla_{\bm{x}}\log\kappa = \hat{\bm{\Phi}}_{{tt_0}}^{\top}\hat{\bm{\Gamma}}_{{tt_0}}^{-1}\left(\bm{y}-\hat{\bm{\Phi}}_{{tt_0}}\bm{x}\right)-\bm{\Pi}(t_0, \bm{0}, t)\bm{x},
\label{GradOfLog}    
\end{align}
and since $\nabla_{\bm{x}}^{2} = \nabla_{\bm{x}}\circ \nabla_{\bm{x}}$,
\begin{align}
\nabla_{\bm{x}}^{2}\log\kappa = -\hat{\bm{\Phi}}_{{tt_0}}^{\top}\hat{\bm{\Gamma}}_{{tt_0}}^{-1}\hat{\bm{\Phi}}_{{tt_0}}-\bm{\Pi}(t_0, \bm{0}, t).
\label{HessOfLog}    
\end{align}
We next use the identity
\begin{align}
\nabla_{\bm{x}}^{2}\log\kappa &= \frac{\nabla_{\bm{x}}^{2}\kappa}{\kappa} - \left(\nabla_{\bm{x}}\log\kappa\right)\left(\nabla_{\bm{x}}\log\kappa\right)^{\top}\nonumber\\
\Rightarrow \nabla_{\bm{x}}^{2}\kappa &= \kappa\bigg\{\left(\nabla_{\bm{x}}\log\kappa\right)\left(\nabla_{\bm{x}}\log\kappa\right)^{\top}+\nabla_{\bm{x}}^{2}\log\kappa\bigg\}.
\label{IdentityHessOfLog}
\end{align}
Substituting \eqref{GradOfLog} and \eqref{HessOfLog} in \eqref{IdentityHessOfLog}, we obtain
\begin{align}
&\nabla_{\bm{x}}^{2}\kappa = \kappa\bigg\{\left(\hat{\bm{\Phi}}_{{tt_0}}^{\top}\hat{\bm{\Gamma}}_{{tt_0}}^{-1}\left(\bm{y}-\hat{\bm{\Phi}}_{{tt_0}}\bm{x}\right)-\bm{\Pi}(t_0, \bm{0}, t)\bm{x}\right)\nonumber\\
&\times\!\left(\hat{\bm{\Phi}}_{{tt_0}}^{\top}\hat{\bm{\Gamma}}_{{tt_0}}^{-1}\left(\bm{y}-\hat{\bm{\Phi}}_{{tt_0}}\bm{x}\right)-\bm{\Pi}(t_0, \bm{0}, t)\bm{x}\right)^{\!\top}\!\!-\hat{\bm{\Phi}}_{{tt_0}}^{\top}\hat{\bm{\Gamma}}_{{tt_0}}^{-1}\hat{\bm{\Phi}}_{{tt_0}}\nonumber\\
&-\bm{\Pi}(t_0, \bm{0}, t)\bigg\},
\label{HessOfKappa}
\end{align}
and thus, the diffusion term in the RHS of \eqref{AdvectionReactionDiffusionPDE} is
\begin{align}
&\langle\bm{B}_t \bm{B}_t^{\top},\nabla_{\bm{x}}^{2}\kappa\rangle\nonumber\\
=&\kappa\langle\bm{B}_t \bm{B}_t^{\top},\hat{\bm{\Phi}}_{{tt_0}}^{\top}\hat{\bm{\Gamma}}_{{tt_0}}^{-1}(\bm{y}-\hat{\bm{\Phi}}_{{tt_0}}\bm{x})(\bm{y}-\hat{\bm{\Phi}}_{{tt_0}}\bm{x})^{\top}\hat{\bm{\Gamma}}_{{tt_0}}^{-1}\hat{\bm{\Phi}}_{{tt_0}}\nonumber\\
&-\bm{\Pi x}(\bm{y}-\hat{\bm{\Phi}}_{{tt_0}}\bm{x})^{\top}\hat{\bm{\Gamma}}_{{tt_0}}^{-1}\hat{\bm{\Phi}}_{{tt_0}} - \hat{\bm{\Phi}}_{{tt_0}}^{\top}\hat{\bm{\Gamma}}_{{tt_0}}^{-1}(\bm{y}-\hat{\bm{\Phi}}_{{tt_0}}\bm{x})\bm{x}^{\top}\bm{\Pi} \nonumber\\
&+\bm{\Pi}\bm{xx}^{\top}\bm{\Pi} -\hat{\bm{\Phi}}_{{tt_0}}^{\top}\hat{\bm{\Gamma}}_{{tt_0}}^{-1}\hat{\bm{\Phi}}_{{tt_0}}-\bm{\Pi}(t_0, \bm{0}, t)\rangle.
\label{FrobeniusInnerProduct}
\end{align}

We note that the advection term in the RHS of \eqref{AdvectionReactionDiffusionPDE} is 
\begin{align}
-\langle\nabla_{\bm{x}},\kappa\:\bm{A}_{t}\bm{x}\rangle = -\kappa\langle \nabla_{\bm{x}}\log\kappa,\bm{A}_t\bm{x}\rangle - \kappa\:\tr\bm{A}_t.
\label{advectiontermsiplified}    
\end{align}
Substituting \eqref{LHSadvectionreactiondiffusionPDE}, \eqref{GradOfLog}, \eqref{FrobeniusInnerProduct} and \eqref{advectiontermsiplified} in \eqref{AdvectionReactionDiffusionPDE}, we arrive at
\begin{align}
&\frac{\dot{c}}{c} - \frac{1}{2}\begin{pmatrix}
\bm{x}\\
\bm{y}
\end{pmatrix}^{\top}\dot{\bm{M}}_{tt_0}\begin{pmatrix}
\bm{x}\\
\bm{y}
\end{pmatrix}\nonumber\\
&=-\langle\hat{\bm{\Phi}}_{{tt_0}}^{\top}\hat{\bm{\Gamma}}_{{tt_0}}^{-1}(\bm{y}-\hat{\bm{\Phi}}_{{tt_0}}\bm{x}),\bm{A}_t\bm{x}\rangle -\tr\bm{A}_t\nonumber\\
&+\langle\bm{B}_t \bm{B}_t^{\top},\hat{\bm{\Phi}}_{{tt_0}}^{\top}\hat{\bm{\Gamma}}_{{tt_0}}^{-1}(\bm{y}-\hat{\bm{\Phi}}_{{tt_0}}\bm{x})(\bm{y}-\hat{\bm{\Phi}}_{{tt_0}}\bm{x})^{\top}\hat{\bm{\Gamma}}_{{tt_0}}^{-1}\hat{\bm{\Phi}}_{{tt_0}}\nonumber\\
&-\bm{\Pi x}(\bm{y}-\hat{\bm{\Phi}}_{{tt_0}}\bm{x})^{\top}\hat{\bm{\Gamma}}_{{tt_0}}^{-1}\hat{\bm{\Phi}}_{{tt_0}} - \hat{\bm{\Phi}}_{{tt_0}}^{\top}\hat{\bm{\Gamma}}_{{tt_0}}^{-1}(\bm{y}-\hat{\bm{\Phi}}_{{tt_0}}\bm{x})\bm{x}^{\top}\bm{\Pi} \nonumber\\
&+\bm{\Pi}\bm{xx}^{\top}\bm{\Pi} -\hat{\bm{\Phi}}_{{tt_0}}^{\top}\hat{\bm{\Gamma}}_{{tt_0}}^{-1}\hat{\bm{\Phi}}_{{tt_0}}-\bm{\Pi}(t_0, \bm{0}, t)\rangle -\frac{1}{2}\bm{x}^{\top}\bm{Q}_{t}\bm{x}\nonumber\\
&= -\theta(t) - \frac{1}{2}\begin{pmatrix}
\bm{x}\\
\bm{y}
\end{pmatrix}^{\top}
\bm{S}_{tt_0}
\begin{pmatrix}
\bm{x}\\
\bm{y}
\end{pmatrix},
\label{TemporalAndSpatialEquality}
\end{align}
where the last equality grouped the spatially independent terms and spatially dependent (quadratic in $\bm{x},\bm{y}$) terms, then used the definition of $\theta(t)$ from \eqref{deftheta}. In particular, the matrix
$$\bm{S}_{tt_0}=\begin{bmatrix}\bm{S}_{11}(t,t_0) & \bm{S}_{12}(t,t_0)\\
\bm{S}_{12}^{\top}(t,t_0) & \bm{S}_{22}(t,t_0)
\end{bmatrix}\in\mathbb{R}^{2n\times 2n}$$
comprises of $n\times n$ blocks
\begin{subequations}
\begin{align}
\bm{S}_{11}(t,t_0) :=& -2\hat{\bm{\Phi}}^{\top}_{tt_0}\hat{\bm{\Gamma}}_{tt_0}^{-1}\hat{\bm{\Phi}}_{tt_0}\hat{\bm{A}}_{t}\nonumber\\
&-\hat{\bm{\Phi}}^{\top}_{tt_0}\hat{\bm{\Gamma}}_{tt_0}^{-1}\hat{\bm{\Phi}}_{tt_0}\hat{\bm{B}}_{t}\hat{\bm{B}}_{t}^{\top}\hat{\bm{\Phi}}^{\top}_{tt_0}\hat{\bm{\Gamma}}_{tt_0}^{-1}\hat{\bm{\Phi}}_{tt_0}  \nonumber\\
&- \bm{\Pi}(t_0, \bm{0}, t)\hat{\bm{B}}_{t}\hat{\bm{B}}_{t}^{\top}\bm{\Pi}(t_0, \bm{0}, t) + \bm{Q}_{t}, \label{S11}\\
\bm{S}_{12}(t,t_0) := &\hat{\bm{\Gamma}}_{tt_0}^{-1}\hat{\bm{\Phi}}_{tt_0}\left(\hat{\bm{A}}_{t} + \hat{\bm{B}}_{t}\hat{\bm{B}}_{t}^{\top}\hat{\bm{\Phi}}^{\top}_{tt_0}\hat{\bm{\Gamma}}_{tt_0}^{-1}\hat{\bm{\Phi}}_{tt_0}\right), \label{S12}\\
\bm{S}_{22}(t,t_0) :=& -\hat{\bm{\Gamma}}_{tt_0}^{-1}\hat{\bm{\Phi}}_{tt_0}\hat{\bm{B}}_{t}\hat{\bm{B}}_{t}^{\top}\hat{\bm{\Phi}}^{\top}_{tt_0}\hat{\bm{\Gamma}}_{tt_0}^{-1}.\label{S22}
\end{align}
\label{BlocksOfS}
\end{subequations}

Since \eqref{TemporalAndSpatialEquality} holds for arbitrary $\bm{x},\bm{y}\in\mathbb{R}^{n}$, equating the $\bm{x},\bm{y}$ independent terms from both sides of \eqref{TemporalAndSpatialEquality}, we conclude that $c$ solves the linear ODE
\begin{align}
\dot{c} = -\theta(t) c.
\label{cODEforLQ}
\end{align}
The solution of \eqref{cODEforLQ} is
\begin{align}
c(t,t_0) = a\exp\left(-\int_{t_0}^{t}\theta(s)\differential s\right),
\label{cUptoConstanta}    
\end{align}
for to-be-determined constant $a$. In Step 2 detailed next, we determine $a$. 

As a by-product of the above calculation, we get $\dot{\bm{M}}_{tt_0}=\bm{S}_{tt_0}$, but this will not be used hereafter.

\noindent{\textbf{Step 2: Determining the constant $a$.}}\\
\noindent 
Substituting \eqref{cUptoConstanta} and \eqref{HalfDistSquaredABQ}-\eqref{defMtt0} in \eqref{kappaExpDist}, letting $t\downarrow t_0$, and invoking the initial condition \eqref{AdvectionReactionDiffusionIC}, we get
\begin{align}
\delta(\bm{x}-\bm{y}) = a\lim_{t\downarrow t_0}\exp\left(-\int_{t_0}^{t}\theta(\tau)\differential\tau-\frac{1}{2}\!\begin{pmatrix}
\bm{x}\\ \bm{y}
\end{pmatrix}^{\!\!\top}\bm{M}_{tt_{0}}\begin{pmatrix}
\bm{x}\\ \bm{y}
\end{pmatrix}\right).
\label{FindingaLQstart}
\end{align}
Integrating both sides of \eqref{FindingaLQstart} w.r.t. $\bm{x}$, and using the shift property of the Dirac delta: $\int_{\mathbb{R}^{n}}\delta(\bm{x}-\bm{y})\cdot 1\cdot\differential\bm{x} = 1$, we obtain 
\begin{align}
1 = a\int_{\mathbb{R}^{n}}\lim_{t\downarrow t_0}\exp\!\left(\!-\!\int_{t_0}^{t}\!\theta(\tau)\differential\tau\!\right)\!\exp\!\left(\!\!-\frac{1}{2}\!\begin{pmatrix}
\bm{x}\\ \bm{y}
\end{pmatrix}^{\!\!\top}\!\!\bm{M}_{tt_{0}}\!\!\begin{pmatrix}
\bm{x}\\ \bm{y}
\end{pmatrix}\!\!\right)\!\differential\bm{x}.
\label{FindingaLQmiddle}
\end{align}
As noted in Sec. \ref{subsec:gettingdist}, $\bm{M}_{tt_0}\succ\bm{0}$, so for any $\bm{x},\bm{y}\in\mathbb{R}^{n}$, the second exponential term in \eqref{FindingaLQmiddle} takes values in $[0,1]$.

Also, under assumption \textbf{A1}, $\theta(\tau)$ in \eqref{deftheta} is bounded for $0\leq t_{0}< t < \infty$, and so is $\exp(-\int_{t_0}^{t}\theta(\tau)\differential\tau)$. Thus, the product of the two exponential terms in \eqref{FindingaLQmiddle} is bounded w.r.t. $\bm{x}$. Therefore, by dominated convergence theorem, we exchange the order of integral and limit in \eqref{FindingaLQmiddle}, to find
\begin{align}
&1= a\lim_{t\downarrow t_0}\exp\!\left(\!-\!\int_{t_0}^{t}\!\theta(\tau)\differential\tau\!\right)\!\int_{\mathbb{R}^{n}}\!\exp\!\left(\!\!-\frac{1}{2}\!\begin{pmatrix}
\bm{x}\\ \bm{y}
\end{pmatrix}^{\!\!\top}\!\!\bm{M}_{tt_{0}}\!\!\begin{pmatrix}
\bm{x}\\ \bm{y}
\end{pmatrix}\!\!\right)\!\differential\bm{x}\nonumber\\
&= a\lim_{t\downarrow t_0}\exp\!\left(\!-\!\int_{t_0}^{t}\!\theta(\tau)\differential\tau\!\right)\!\exp\!\left(\!-\frac{1}{2}\bm{y}^{\top}\bm{M}_{22}(t,t_0)\bm{y}\!\right)\nonumber\\
&\quad\times\int_{\mathbb{R}^{n}}\!\!\exp\!\left(\!-\frac{1}{2}\bm{x}^{\top}\bm{M}_{11}(t,t_0)\bm{x}+\langle-\bm{M}_{12}(t,t_0)\bm{y},\bm{x}\rangle\vphantom{\frac{1}{2}}\right)\differential\bm{x}\nonumber\\
&= a\left(2\pi\right)^{\frac{n}{2}}\lim_{t\downarrow t_0}\exp\!\left(\!-\!\int_{t_0}^{t}\!\theta(\tau)\differential\tau\!\right)\!\left(\det\bm{M}_{11}(t,t_0)\right)^{-\frac{1}{2}}\nonumber\\
&\exp\!\left(\!-\frac{1}{2}\bm{y}^{\top}\{\bm{M}_{22}(t,t_0)-\bm{M}_{12}^{\top}(t,t_0)\bm{M}_{11}^{-1}(t,t_0)\bm{M}_{12}(t,t_0)\}\bm{y}\vphantom{\frac{1}{2}}\!\right),
\label{FindingaLQEND}
\end{align}
wherein $\bm{M}_{11},\bm{M}_{12},\bm{M}_{22}$ are the respective $n\times n$ blocks of \eqref{defMtt0}, and the last step follows from Lemma \ref{Lemma:GaussianIntegral}. Invoking Lemma \ref{Lemma:GaussianIntegral} here is in turn made possible by the fact that $\bm{M}_{11}$ being a sum of two positive definite matrices (see \eqref{defMtt0}), is positive definite.

Now the idea is to evaluate the limits in \eqref{FindingaLQEND}. For the exponential of negative quadratic term, using \eqref{defMtt0}, we find 
{\small{\begin{align}
&\bm{M}_{22}(t,t_0)-\bm{M}_{12}^{\top}(t,t_0)\bm{M}_{11}^{-1}(t,t_0)\bm{M}_{12}(t,t_0)\nonumber\\
=& \hat{\bm{\Gamma}}_{tt_0}^{-1}-\hat{\bm{\Gamma}}^{-1}_{tt_0}\hat{\bm{\Phi}}_{tt_0} \!\left(\hat{\bm{\Phi}}_{tt_0}^{\top}\hat{\bm{\Gamma}}_{tt_0}^{-1}\hat{\bm{\Phi}}_{tt_0}+\bm{\Pi}(t_0, \bm{0}, t)\right)^{\!-1}\hat{\bm{\Phi}}_{tt_0}^{\top}\hat{\bm{\Gamma}}^{-1}\!\nonumber\\
=& \hat{\bm{\Gamma}}_{tt_0}^{-1} \notag\\
&\!-\left(\hat{\bm{\Phi}}_{tt_0}^{-1}\hat{\bm{\Gamma}}_{tt_0}\right)^{\!-1}\!\!\left(\hat{\bm{\Phi}}_{tt_0}^{\top}\hat{\bm{\Gamma}}_{tt_0}^{-1}\hat{\bm{\Phi}}_{tt_0}+\bm{\Pi}(t_0, \bm{0}, t)\right)^{\!-1}\!\left(\hat{\bm{\Gamma}}_{tt_0}\hat{\bm{\Phi}}_{tt_0}^{-\top}\right)^{-1}\nonumber\\
=&\hat{\bm{\Gamma}}_{tt_0}^{-1} - \left(\hat{\bm{\Gamma}}_{tt_0} + \hat{\bm{\Gamma}}_{tt_0}\hat{\bm{\Phi}}_{tt_0}^{-\top}\bm{\Pi}(t_0, \bm{0}, t)\hat{\bm{\Phi}}_{tt_0}^{-1}\hat{\bm{\Gamma}}_{tt_0}\right)^{-1},
\label{SchurComplementWithinExp}    
\end{align}
}}
thanks to the invertibility of $\hat{\bm{\Phi}}_{tt_0},\hat{\bm{\Gamma}}_{tt_0}$. Using the Woodbury identity \cite{hager1989updating}, \cite[p. 19]{horn2012matrix},
\begin{align}
&\left(\hat{\bm{\Gamma}}_{tt_0} + \hat{\bm{\Gamma}}_{tt_0}\hat{\bm{\Phi}}_{tt_0}^{-\top}\bm{\Pi}(t_0, \bm{0}, t)\hat{\bm{\Phi}}_{tt_0}^{-1}\hat{\bm{\Gamma}}_{tt_0}\right)^{-1}\nonumber\\
=& \hat{\bm{\Gamma}}_{tt_0}^{-1} - \hat{\bm{\Phi}}_{tt_0}^{-\top}\left(\bm{\Pi}^{-1}(t_0, \bm{0}, t) + \hat{\bm{\Phi}}_{tt_0}^{-1}\hat{\bm{\Gamma}}_{tt_0}\hat{\bm{\Phi}}_{tt_0}^{-\top} \right)^{-1}\hat{\bm{\Phi}}_{tt_0}^{-1},
\label{ApplyingWoodbury}
\end{align}
which upon substituting in  \eqref{SchurComplementWithinExp}, gives
\begin{align}
&\bm{M}_{22}(t,t_0)-\bm{M}_{12}^{\top}(t,t_0)\bm{M}_{11}^{-1}(t,t_0)\bm{M}_{12}(t,t_0)\nonumber\\
=&\hat{\bm{\Phi}}_{tt_0}^{-\top}\left(\bm{\Pi}^{-1}(t_0, \bm{0}, t) + \hat{\bm{\Phi}}_{tt_0}^{-1}\hat{\bm{\Gamma}}_{tt_0}\hat{\bm{\Phi}}_{tt_0}^{-\top} \right)^{-1}\hat{\bm{\Phi}}_{tt_0}^{-1}.
\label{SchurComplementSimplified}
\end{align}
Since $\hat{\bm{\Phi}}_{t_0 t_0}=\bm{I}$, $\hat{\bm{\Gamma}}_{t_0 t_0}=\bm{0}$, $\bm{\Pi}(t_0, \bm{0}, t_0)=\bm{0}$, we then have
{\small{
\begin{align}
&\lim_{t\downarrow t_0}\exp\left(-\frac{1}{2}\bm{y}^{\top}\{\bm{M}_{22}(t,t_0)-\bm{M}_{12}^{\top}(t,t_0)\bm{M}_{11}^{-1}(t,t_0)\bm{M}_{12}(t,t_0)\}\bm{y}\right)\nonumber\\
&=\exp\!\!\left(\!\!-\frac{1}{2}\bm{y}^{\top}\underbrace{\lim_{t\downarrow t_0}\{\hat{\bm{\Phi}}_{tt_0}^{-\top}\!\!\left(\!\bm{\Pi}^{-1}(t_0, \bm{0}, t) + \hat{\bm{\Phi}}_{tt_0}^{-1}\hat{\bm{\Gamma}}_{tt_0}\hat{\bm{\Phi}}_{tt_0}^{-\top} \right)^{\!-1}\!\hat{\bm{\Phi}}_{tt_0}^{-1}\}}_{=\bm{0}}\bm{y}\!\!\right)\nonumber\\
&= 1.
\label{LimitForExpTerm}
\end{align}
}}
Therefore, \eqref{FindingaLQEND} simplifies to 
\begin{align}
1 = a\left(2\pi\right)^{\frac{n}{2}}\lim_{t\downarrow t_0}\exp\!\left(\!-\!\int_{t_0}^{t}\!\theta(\tau)\differential\tau\!\right)\!\left(\det\bm{M}_{11}(t,t_0)\right)^{-\frac{1}{2}}.
\label{aEquality}    
\end{align}
Notice that for the limit in \eqref{aEquality} to be defined, the limit cannot be further distributed to the exponential and determinant factors. 

Since the limit of reciprocal equals to the reciprocal of the limit, and by Jacobi identity: $\exp\tr(\cdot)=\det\exp(\cdot)$, we re-write \eqref{aEquality} as \eqref{defa}.

\noindent{\textbf{Step 3: Putting everything together.}}\\
Substituting \eqref{defa} in \eqref{cUptoConstanta}, and then substituting the resulting expression for $c(t,t_0)$ in \eqref{kappaintermedform}, the statement follows.\qed


\subsection{Proof of Corollary \ref{Corollary:SplCaseOfLQkernel}}\label{AppProofOfCorollarySplCase}
Under the stated conditions, we proved in Proposition \ref{prop:SpecialCaseAzeroBidentity} that the term $\frac{1}{2}\!\begin{pmatrix}
\bm{x}\\ \bm{y}
\end{pmatrix}^{\!\!\top}\!\!\bm{M}_{tt_{0}}\!\!\begin{pmatrix}
\bm{x}\\ \bm{y}
\end{pmatrix}$ in \eqref{KernelInTermsOfM} reduces to \eqref{halfdistsquaredSICONkernel}. What remains is to show that the term $\exp\left(-\int_{t_0}^{t}\theta(s)\differential s\right)$ in \eqref{KernelInTermsOfM} reduces to $\displaystyle\prod_{i=1}^{n}\dfrac{1}{\left(\sinh(\omega_i (t-t_0))\right)^{1/2}}$ as found in \eqref{cwitha}, and to compute the pre-factor $a$.

To this end, notice that \eqref{deftheta} in this case specializes to
\begin{align}
\theta(s) &= \tr\left(\hat{\bm{\Phi}}_{{s t_0}}^{\top}\hat{\bm{\Gamma}}_{{s t_0}}^{-1}\hat{\bm{\Phi}}_{{s t_0}} + \bm{\Pi}(t_0, \bm{0}, s)\right)\nonumber\\
&= \tr\left(\sqrt{\bm{D}}\coth\left(2\sqrt{\bm{D}}(s-t_0)\right)\right)\nonumber\\
&= \sum_{i=1}^{n}\sqrt{D_{ii}}\coth\left(\sqrt{D_{ii}}(s-t_0)\right),
\label{thetasplcase}
\end{align}
where the second equality follows from \eqref{Pisplcase}, \eqref{SpecialCasePhihat} and \eqref{SpecialCaseInvGammahat}. Then
\begin{align}
&\exp\left(-\int_{t_0}^{t}\theta(s)\differential s\right)\nonumber\\
=&\lim_{\varepsilon\downarrow 0}\exp\left(-\sum_{i=1}^{n}\sqrt{D}_{ii}\int_{\varepsilon}^{t-t_0}\coth\left(2\sqrt{D}_{ii}\tau\right)\differential\tau\right)\nonumber\\
=&\exp\left(-\sum_{i=1}^{n}\frac{1}{2}\log\sinh\left(2\sqrt{D}_{ii}(t-t_0)\right)\right)\nonumber\\
=&\displaystyle\prod_{i=1}^{n}\dfrac{1}{\left(\sinh(\omega_i (t-t_0))\right)^{1/2}}.
\label{expMinusIntegralTheta}
\end{align}
In above, the first equality used \eqref{thetasplcase} and a change-of-variable $s-t_0 \mapsto \tau$. The second equality follows from integration and the limit. The last equality used $\omega_i := 2\sqrt{D}_{ii}$. 

To compute the pre-factor $a$ for this case, we note from \eqref{Nsplcase} that $\det\bm{M}_{11}^{1/2}(t,t_0)=\bm{D}^{1/4}\left(\coth\left(2\sqrt{\bm{D}}(t-t_0)\right)\right)^{1/2}$. Then using the Jacobi identity $\det\exp(\cdot) = \exp\tr(\cdot)$, we have
\begin{align}
a &= (2\pi)^{-n/2}\lim_{t\downarrow t_0}\bigg\{\det\bm{M}_{11}^{1/2}(t,t_0)\left(\exp\int_{t_0}^{t}\theta(\tau)\differential\tau\right)\bigg\}\nonumber\\
&= (2\pi)^{-n/2}\bm{D}^{1/4}\lim_{t\downarrow t_0}\bigg\{\left(\coth\left(2\sqrt{\bm{D}}(t-t_0)\right)\right)^{1/2}\nonumber\\
&\qquad\qquad\qquad\qquad\quad\;\times\left(\sinh\left(2\sqrt{\bm{D}}(t-t_0)\right)\right)^{1/2}\bigg\}\nonumber\\
&=(2\pi)^{-n/2}\bm{D}^{1/4}\lim_{t\downarrow t_0}\left(\cosh\left(2\sqrt{\bm{D}}(t-t_0)\right)\right)^{1/2}\nonumber\\
&=(2\pi)^{-n/2}\bm{D}^{1/4},
\label{aSplCase}
\end{align}
where the second equality is due to \eqref{expMinusIntegralTheta}. 

With the $a$ given by \eqref{aSplCase}, $\exp(-\int_{t_0}^{t}\theta(s)\differential s)$ given by \eqref{expMinusIntegralTheta}, and the specialization of $\bm{M}_{tt_0}$ as in Proposition \ref{prop:SpecialCaseAzeroBidentity}, the kernel \eqref{KernelInTermsOfM} specializes to \cite[Eq. (43)]{teter2024weyl} or that in \cite[Eq. (A.22)]{teter2024schr}.\qed

\section*{References}
\bibliographystyle{IEEEtran}
\bibliography{References.bib}

\balance


\end{document}